\documentclass[11pt,a4paper]{article}

\usepackage[T1]{fontenc}
\usepackage[utf8]{inputenc}
\usepackage[english]{babel}
\usepackage{lmodern}
\usepackage{microtype}
\usepackage{amsmath,amssymb,amsthm,mathtools}
\usepackage{mathrsfs}
\usepackage[margin=1.15in]{geometry}
\usepackage{float}
\usepackage[hidelinks]{hyperref}
\usepackage{enumitem}
\usepackage{tikz}
\newtheorem{theorem}{Theorem}

\newtheorem{proposition}{Proposition}[section]
\newtheorem{lemma}[proposition]{Lemma}
\newtheorem{corollary}[proposition]{Corollary}
\theoremstyle{definition}
\newtheorem{definition}[proposition]{Definition}

\theoremstyle{remark}
\newtheorem{remark}[proposition]{Remark}

\theoremstyle{plain}
\newtheorem{claim}{Claim}

\AtEndEnvironment{proof}{\setcounter{claim}{0}}
\newenvironment{theoremproof}[1]
  {\begin{proof}[\textbf{\textit{Proof of Theorem~\ref{#1}}}]}
  {\end{proof}}

\renewcommand{\varnothing}{\emptyset}

\newcommand{\Cantor}{2^\omega}
\newcommand{\Baire}{\omega^\omega}
\newcommand{\Fin}{[\omega]^{<\omega}}
\newcommand{\Inf}{[\omega]^\omega}
\newcommand{\FU}{\operatorname{FU}}
\newcommand{\Pow}{\mathcal P}
\newcommand{\restr}{\mathbin{\upharpoonright}}
\newcommand{\non}{\operatorname{non}}
\newcommand{\wlambda}{\mathrm{w}\lambda}
\newcommand{\CH}{\mathsf{CH}}

\DeclareMathOperator{\Con}{Con}
\newcommand{\ZFC}{\mathsf{ZFC}}
\newcommand{\MA}{\mathsf{MA}}
\newcommand{\OCA}{\mathsf{OCA}_{\mathrm T}}

\newcommand{\A}{\mathcal A}
\newcommand{\B}{\mathcal B}
\newcommand{\C}{\mathcal C}
\newcommand{\D}{\mathcal D}
\newcommand{\F}{\mathcal F}
\newcommand{\cl}{\operatorname{cl}}
\newcommand{\Fn}{\operatorname{Fn}}
\newcommand{\eqstar}{\mathrel{=^*}}
\def\cA{{\mathcal{A}}}
\def\cB{{\mathcal{B}}}
\def\cD{{\mathcal{D}}}

\title{Weak \texorpdfstring{$\lambda$}{lambda}-sets and almost-normal
Isbell--Mrówka spaces}
\author{C\'esar Ismael Corral Rojas\thanks{Departamento de Matem\'aticas,
Facultad de Ciencias, Universidad Nacional Aut\'onoma de M\'exico, Mexico City,
Mexico. Email: \texttt{cicorral@ciencias.unam.mx}.}
\and Vinicius de Oliveira Rodrigues\thanks{Departamento de Matem\'atica,
Instituto de Matem\'atica, Estat\'istica e Ci\^encia da Computa\c{c}\~ao,
Universidade de S\~ao Paulo,
S\~ao Paulo, Brazil. Email: \texttt{vinior@ime.usp.br}.}
\and Victor dos Santos Ronchim\thanks{Departamento de Bioprocessos e
Biotecnologia, Faculdade de Ci\^encias Agron\^omicas, Universidade Estadual
Paulista, Botucatu, Brazil. Email: \texttt{victor.ronchim@unesp.br}.}}
\date{}

\hypersetup{
  pdftitle={Weak lambda-sets and almost-normal Isbell--Mrówka spaces},
  pdfsubject={Set-theoretic topology},
  pdfkeywords={weak lambda-set, almost-normality, almost disjoint family,
    Hindman's theorem, Open Coloring Axiom}
}

\begin{document}

\begingroup
\makeatletter
\let\@fnsymbol\@arabic
\makeatother
\maketitle
\endgroup

\begin{abstract}
We determine the least cardinality of a set of reals that is not a weak
$\lambda$-set, proving that it is $\max\{\mathfrak b,\mathfrak s\}$, where
$\mathfrak b$ and $\mathfrak s$ are the bounding and splitting numbers.
We prove that if $\mathfrak s<\mathfrak b$, then every first countable weak $\lambda$-space
is a $\lambda$-space.  We also show, relative to the consistency of $\ZFC$,
that the existence of a weak $\lambda$-set that is not a $\lambda$-set is
independent of $\ZFC+\mathfrak b=\mathfrak s$.
In particular, we show that in the model obtained by adding $\omega_2$
Cohen reals over a model of $\CH$, every weak $\lambda$-set is a $\lambda$-set.
Finally, we construct an almost-normal
MAD family under $\CH$ and show that no such family exists
under $\OCA+\mathfrak{ap}>\omega_1$.
Thus, the existence of an almost-normal MAD family is independent of $\ZFC$.
\end{abstract}



\section{Introduction}

The study of special subsets of the real line is a classical topic at the
intersection of topology and set theory.
Such sets are defined by properties concerning descriptive complexity, Baire category and measure theoretic notions.  For example, a $\lambda$-set is a set in which every
countable subset is a relative $G_\delta$, a $\sigma$-set is one in which every
$G_\delta$ subset is a relative $F_\sigma$, and a $Q$-set is one in which every subset
is a relative $G_\delta$.  The possible
cardinalities of these sets and the relationships between them are closely tied to
cardinal characteristics of the continuum and to additional set-theoretic
assumptions.  We refer to \cite{MillerSurvey} for an overview of this field of research.

These classes are also deeply connected to the study of separation properties of
topological spaces, such as normality. This connection led to the introduction of weak
$\lambda$-sets in \cite{RodriguesRonchimSzeptycki}, a weakening of $\lambda$-sets, while studying separation properties of Isbell--Mrówka spaces.

\begin{definition}
\label{def:weak-lambda-space}
A topological space $X$ is a \emph{weak $\lambda$-space} if, whenever
$A,B\subseteq X$ are disjoint and countable, there is a set $H\subseteq X$
which is both $F_\sigma$ and $G_\delta$ in $X$ such that
$A\subseteq H$ and $B\subseteq X\setminus H$.  A separable metrizable weak
$\lambda$-space is called a \emph{weak $\lambda$-set}.
\end{definition}

Every $\lambda$-set is a weak $\lambda$-set.  Rodrigues, Ronchim and
Szeptycki asked whether the converse can consistently fail
\cite[Question~4.6]{RodriguesRonchimSzeptycki}.

An almost disjoint family is a collection $\A\subseteq[\omega]^\omega$
whose distinct members have finite intersection.  Given such a family, the
Isbell--Mrówka space $\Psi(\A)$ is obtained by adjoining to the discrete
space $\omega$ one point for each member of $\A$.  A basic neighborhood of
$a\in\A$ is of the form $\{a\}\cup(a\setminus F)$ for some finite $F\subseteq a$.  The
space $\Psi(\A)$ is Hausdorff, separable, scattered, zero-dimensional,
first countable, locally compact and regular.  Its topological properties
are often related to combinatorial properties of the underlying almost disjoint family $\A$: for instance, $\Psi(\A)$ is
pseudocompact if and only if $\A$ is maximal, and $\Psi(\mathcal A)$ is metrizable if and only
if $\A$ is countable.  An infinite almost disjoint family $\A$ is maximal
(a MAD family) if every infinite subset of $\omega$ meets some member of
$\A$ infinitely.
For background on Isbell--Mrówka spaces, see
\cite{HernandezHrusakMrowkaIsbell,HrusakAlmostDisjointTopology}.

Since $\A$ is closed and discrete, $\Psi(\A)$ is not countably compact
when $\A$ is infinite \cite{HernandezHrusakMrowkaIsbell,HrusakAlmostDisjointTopology}.
A normal pseudocompact space is countably compact, so $\Psi(\A)$ is not
normal when $\A$ is MAD. Also, Jones' Lemma shows that $\Psi(\A)$ is not
normal whenever $|\A|=\mathfrak c$.
Moreover, in $\ZFC$, there exist almost disjoint families of cardinality
$\omega_1$ which are not normal, e.g., Luzin families
(see, e.g., \cite[Sections~3.1--3.2]{HrusakAlmostDisjointTopology}).  Nevertheless, it is consistent with $\ZFC$
that there are uncountable almost disjoint families $\A$ for which
$\Psi(\A)$ is normal \cite{TallNormalMoore}.

Several weakenings of normality for $\Psi$-spaces have been studied,
such as almost-normality, quasi-normality, partly-normality, and mild
normality \cite{GarciaBalanSzeptycki}.
In this paper, we focus on almost-normality.  A space is almost-normal if any two disjoint
closed sets, one of which is regular closed, have disjoint open
neighborhoods.

For $X\subseteq2^\omega$, let $\A_X$ consist of the
branches $\{x\restr n:n\in\omega\}$, for $x\in X$.  The space
$\Psi(\A_X)$ is normal if and only if $X$ is a $Q$-set, almost-normal if
and only if $X$ is a $\sigma$-set, and pseudonormal if and only if $X$ is
a $\lambda$-set \cite[Corollary~2.3, Theorem~2.9 and Proposition~3.3]{RodriguesRonchimSzeptycki}.

Another separation property related to Isbell--Mrówka spaces is strong
$\aleph_0$-separation.  We say that $\A$ is strongly $\aleph_0$-separated
if any two disjoint countable closed subsets of $\Psi(\A)$ can be separated
by clopen sets.  For $\A_X$ as above, this is equivalent to $X$ being a
weak $\lambda$-set \cite[Proposition~4.3]{RodriguesRonchimSzeptycki}.

For a class of spaces $\mathcal P$, write $\non(\mathcal P)$ for the least
cardinality of a separable metrizable space outside $\mathcal P$.  Let
$\mathfrak b$ and $\mathfrak s$ denote the bounding and splitting numbers,
respectively (for their definitions, see Section~\ref{sec:preliminaries}).  It is well-known that
$\non(\lambda)=\non(\sigma)=\mathfrak b$
\cite[Theorem~9.1]{vanDouwen}.  We also write
$\mathfrak q=\non(Q)$ for the least cardinality of a set of reals that is
not a $Q$-set \cite[p.~13]{BrendleDow}.
The least cardinality of a set of reals that is not a weak $\lambda$-set, that is, $\non(\wlambda)$, was left open \cite{RodriguesRonchimSzeptycki}.
In this paper, we settle this question by proving that $\non(\wlambda)=\max\{\mathfrak b,\mathfrak s\}$.
\begin{theorem}
\label{thm:intro-non-weak-lambda}
    $\non(\wlambda)=\max\{\mathfrak b,\mathfrak s\}$.
\end{theorem}
The following diagram summarizes these relationships
for families of branches.
\begin{figure}[H]
\centering
\resizebox{\linewidth}{!}{%
\begin{tikzpicture}[
  every node/.style={font=\small},
  implication/.style={double, double distance=0.65pt, -stealth,
    shorten <=2pt, shorten >=2pt},
  correspondence/.style={<->, shorten <=3pt, shorten >=3pt}
]
  \node[anchor=west] at (0,1.9) {$\non(\mathcal P)$:};
  \node[anchor=west] at (0,0.7) {$X\subseteq2^\omega$:};
  \node[anchor=west] at (0,-1.0) {$\Psi(\A_X)$:};

  \node (cardw) at (3.25,1.9) {$\max\{\mathfrak b,\mathfrak s\}$};
  \node (cardl) at (6.6,1.9) {$\mathfrak b$};
  \node (cards) at (9.65,1.9) {$\mathfrak b$};
  \node (cardq) at (12.55,1.9) {$\mathfrak q$};
  \node at (4.95,1.9) {$\geq$};
  \node at (8.12,1.9) {$=$};
  \node at (11.10,1.9) {$\geq$};

  \node (weak) at (3.25,0.7) {weak $\lambda$-set};
  \node (lambda) at (6.6,0.7) {$\lambda$-set};
  \node (sigma) at (9.65,0.7) {$\sigma$-set};
  \node (qset) at (12.55,0.7) {$Q$-set};

  \node[align=center] (separated) at (3.25,-1.0)
    {strongly $\aleph_0$-\\separated};
  \node (pseudo) at (6.6,-1.0) {pseudonormal};
  \node (almost) at (9.65,-1.0) {almost-normal};
  \node (normal) at (12.55,-1.0) {normal};

  \draw[implication] (lambda.west) -- (weak.east);
  \draw[implication] (sigma.west) -- (lambda.east);
  \draw[implication] (qset.west) -- (sigma.east);
  \draw[implication] (pseudo.west) -- (separated.east);
  \draw[implication] (almost.west) -- (pseudo.east);
  \draw[implication] (normal.west) -- (almost.east);

  \draw[correspondence] (weak.south) -- (separated.north);
  \draw[correspondence] (lambda.south) -- (pseudo.north);
  \draw[correspondence] (sigma.south) -- (almost.north);
  \draw[correspondence] (qset.south) -- (normal.north);
\end{tikzpicture}%
}
\caption{Special sets of reals and separation properties of their families
of branches.  The top row gives the least cardinality of a set of reals
outside each class.  Horizontal arrows denote implications and vertical arrows
denote equivalences.}
\label{fig:special-sets-separation}
\end{figure}

Rodrigues, Ronchim and Szeptycki asked whether there is, at least
consistently, a weak $\lambda$-set which is not a $\lambda$-set
\cite[Question~4.6]{RodriguesRonchimSzeptycki}.  The preceding theorem,
together with $\non(\lambda)=\mathfrak b$, gives a positive answer when
$\mathfrak b<\mathfrak s$.
Since each of the inequalities $\mathfrak b<\mathfrak s$ and
$\mathfrak s<\mathfrak b$ is consistent relative to $\ZFC$, the
observation above and the following theorem show that the existence of a
weak $\lambda$-set which is not a $\lambda$-set is independent of $\ZFC$,
relative to its consistency.

\begin{theorem}
\label{thm:intro-s-b-weak-lambda}
If $\mathfrak s<\mathfrak b$, then every weak $\lambda$-set is a
$\lambda$-set.
\end{theorem}

The remaining question is whether there exists a weak $\lambda$-set which is not a $\lambda$-set if $\mathfrak s=\mathfrak b$.
In
the model obtained by adding $\omega_2$ Cohen reals over a model of $\CH$,
we have $\mathfrak s=\mathfrak b=\aleph_1$ \cite[Section~11.3]{Blass}.
Let $\mathbb C_{\omega_2}$ denote the forcing for adding $\omega_2$
Cohen reals.  We prove the following.

\begin{theorem}
\label{thm:intro-cohen-weak-lambda}
Assume $\CH$.  Then $\mathbb C_{\omega_2}\Vdash
\text{every weak $\lambda$-set is a $\lambda$-set}$.
\end{theorem}

On the other hand, we show that if $\mathfrak p=\mathfrak c$, then there is
a weak $\lambda$-set which is not a $\lambda$-set, where $\mathfrak p$ is
the pseudointersection number (for its definition, see
Section~\ref{sec:preliminaries}).
This hypothesis implies $\mathfrak s=\mathfrak b=\mathfrak c$.
Under this assumption, we obtain a positive answer to the question above.
This result was also obtained independently by Osipov and Lipin
\cite[Theorem~1.1]{osipovweaklambda}.

\begin{theorem}
\label{thm:intro-p-c-weak-lambda}
If $\mathfrak p=\mathfrak c$, then there is a weak $\lambda$-set which is
not a $\lambda$-set.
\end{theorem}

Together, these results show that, relative to the consistency of $\ZFC$,
the existence of a weak $\lambda$-set which is not a $\lambda$-set is
independent of $\ZFC+\mathfrak s=\mathfrak b$.

We now turn to almost-normal MAD families.  As noted above, a MAD family
cannot be normal.  Garc\'ia-Balan and Szeptycki asked whether there is an
almost-normal MAD family \cite[Question~4.3]{GarciaBalanSzeptycki}.
The second and third authors used forcing to construct an almost-normal almost
disjoint family which is not normal \cite{RodriguesRonchimAlmostNormal}.
The first author showed that $\CH$ implies the existence of an almost-normal almost
disjoint family which is not normal, and that the Proper Forcing Axiom
implies that no MAD family is almost-normal \cite{CorralPublished}.
An earlier version of the first author's work claimed the existence of an
almost-normal MAD family under $\CH$, but the argument contained a gap
\cite[Theorem~2.7]{CorralPreprint}.  In
Section~\ref{sec:almost-normal-mad-ch}, we
obtain the following result.

\begin{theorem}
\label{thm:intro-ch-almost-normal-mad}
Assuming $\CH$, there exists an almost-normal MAD family.
\end{theorem}

For the opposite direction, we use Todor\v cevi\'c's Open
Coloring Axiom $\OCA$ and the cardinal characteristic of the continuum $\mathfrak{ap}$ \cite[p.~14]{BrendleDow}.  In
Section~\ref{sec:oca-obstruction}, we recall the definition of $\mathfrak{ap}$ and strengthen
the first author's result by proving the following.

\begin{theorem}
\label{thm:intro-oca-ap-no-almost-normal}
Assume $\OCA+\mathfrak{ap}>\omega_1$.  Then almost-normal MAD families do
not exist.
\end{theorem}

Together, the last two theorems stated above show that the existence of an
almost-normal MAD family is independent of $\ZFC$.

\medskip
The paper is organized as follows.
\begin{itemize}[leftmargin=\parindent]
\item Section~\ref{sec:preliminaries}: we introduce the notation and
preliminary results used throughout the paper.

\item Section~\ref{sec:weak-lambda}: we prove
Theorem~\ref{thm:intro-non-weak-lambda}, proving that $\non(\wlambda)=\max\{\mathfrak b,\mathfrak s\}$,
and obtain as a corollary a consistent example of a weak $\lambda$-set which is not a $\lambda$-set under
$\mathfrak b<\mathfrak s$.

\item Section~\ref{sec:direct-ch-weak-lambda}: we prove
Theorem~\ref{thm:intro-p-c-weak-lambda}, constructing a weak $\lambda$-set
which is not a $\lambda$-set under $\mathfrak p=\mathfrak c$.
This result was also obtained independently by Osipov and Lipin
\cite[Theorem~1.1]{osipovweaklambda} using similar techniques before the announcement of the present work.

\item Section~\ref{sec:splitting-bounding-obstruction}: we prove a
strengthening of Theorem~\ref{thm:intro-s-b-weak-lambda}, showing that
$\mathfrak s<\mathfrak b$ implies that every first countable weak
$\lambda$-space is a $\lambda$-space.
This shows that the existence of a weak $\lambda$-set which is not a $\lambda$-set is independent of $\ZFC$.

\item Section~\ref{sec:cohen-model}: we prove
Theorem~\ref{thm:intro-cohen-weak-lambda}, showing that every weak
$\lambda$-set is a $\lambda$-set in the Cohen model.
Together with Theorem~\ref{thm:intro-p-c-weak-lambda}, this yields that the existence of a weak $\lambda$-set which is not a $\lambda$-set is independent of $\ZFC+\mathfrak s=\mathfrak b$.

\item Section~\ref{sec:almost-normal-mad-ch}: we prove
Theorem~\ref{thm:intro-ch-almost-normal-mad}, constructing an
almost-normal MAD family under $\CH$.

\item Section~\ref{sec:oca-obstruction}: we prove
Theorem~\ref{thm:intro-oca-ap-no-almost-normal} showing that no
almost-normal MAD family exists under $\OCA+\mathfrak{ap}>\omega_1$ as a consequence of known results from the literature.
Together with Theorem~\ref{thm:intro-ch-almost-normal-mad}, this implies that the existence of an almost-normal MAD family is independent of $\ZFC$.
\end{itemize}

\section{Preliminaries}
\label{sec:preliminaries}

Our terminology is mainly standard: given two sets $a,b$, we write $a\subseteq^*b$ when $a\setminus b$ is finite and $a\eqstar b$ when $a\mathbin{\triangle}b$ is finite, where $a\mathbin{\triangle}b$ denotes the symmetric difference of $a$ and $b$.  We denote by $\Pow(X)$ the power set of $X$, by $[X]^\kappa$ the set of all subsets of $X$ with cardinality $\kappa$, and $[X]^{<\kappa}=\bigcup_{\eta<\kappa}[X]^\eta$.

For an almost disjoint family over a countably infinite set $N$, we use
$\Psi(\A)$ for the space defined in the introduction, with $N$ in place
of $\omega$.  When forming $N\cup\A$, we assume $N\cap[N]^\omega=\varnothing$.

We will use the following shorthand.

\begin{definition}
    We say that $\A$ satisfies the topological property $P$ if the space $\Psi(\A)$ satisfies $P$.
\end{definition}

Given $\mathscr L, \mathscr R\subseteq\mathscr P(\omega)$, we write $\mathscr L\perp\mathscr R$ if
$|L\cap R|<\omega$ for every $L\in\mathscr L$ and $R\in\mathscr R$.
In this case, we say that $\mathscr L$ and $\mathscr R$ are \emph{orthogonal}.

\begin{definition}
Let $\mathscr L,\mathscr R\subseteq\mathcal P(\omega)$.
A separator for $\mathscr L$ and $\mathscr R$ is a set $D\subseteq\omega$ such that
\[
 L\subseteq^*D\quad(\forall L\in\mathscr L)
 \qquad \text{and}\qquad
 |R\cap D|<\omega\quad(\forall R\in\mathscr R).
\]
If $D$ is a separator for $\mathscr L$ and $\mathscr R$, we say that $D$ separates $\mathscr L$ from $\mathscr R$.
\end{definition}

Notice that if $D$ separates $\mathscr L$ from $\mathscr R$, then $\mathscr L\perp\mathscr R$.
Note that this concept is not symmetric.
However, if $D$ separates $\mathscr L$ from $\mathscr R$, then its complement separates $\mathscr R$ from $\mathscr L$.

It is well known that countable orthogonal families can be separated.
We include its proof for the
convenience of the reader.

\begin{lemma}
\label{lem:countable-separation}
Let $\mathscr L,\mathscr R\subseteq\mathcal P(\omega)$ be countable
collections such that $\mathscr L\perp \mathscr R$, and let $C\in\mathscr L$.
Then there exists a separator $D\supseteq C$ for $\mathscr L$ and $\mathscr R$.
\end{lemma}

\begin{proof}
If $\mathscr R=\varnothing$, take $D=\omega$.  Otherwise, enumerate both
collections with repetitions if needed as $\mathscr L=\{L_n:n<\omega\}$ and $\mathscr R=\{R_n:n<\omega\}$.
Let
\[
\label{eq:diagonal-separator}
 D=C\cup\bigcup_{n<\omega}
 \left(L_n\setminus\bigcup_{j\leq n}R_j\right).
\]
For every $n\in \omega$, $L_n\setminus D\subseteq\bigcup_{j\leq n}(L_n\cap R_j)$,
so $L_n\subseteq^*D$.
On the other hand, for each $m\in \omega$,
\[
 R_m\cap D\subseteq R_m\cap C\cup
\bigcup_{n<m}(R_m\cap L_n).
\]
The right-hand side is finite, so $|R_m\cap D|<\omega$.  Thus $D$ is a separator for $\mathscr L$ and $\mathscr R$.
\end{proof}

A separator that separates an almost disjoint family $\A$ into two orthogonal subfamilies is called a \emph{partitioner} for $\A$.

\begin{definition}
Let $\A$ be an almost disjoint family.
\begin{enumerate}[label=(\arabic*)]
\item A set $D\subseteq\omega$ is a \emph{partitioner} for $\A$ if it separates some subfamily $\B\subseteq\A$ from $\A\setminus\B$.
\item For $C\subseteq\omega$, let $\A_C=\{a\in\A:|a\cap C|=\omega\}$.
\end{enumerate}
\end{definition}

The following proposition reformulates
\cite[Theorem~3.6]{RodriguesRonchimAlmostNormal} and \cite[Proposition~2.3]{CorralPublished}.
Recall from \cite{GarciaBalanSzeptycki} that the
regular closed subsets of $\Psi(\A)$ are precisely the sets
$\cl_{\Psi(\A)}C=C\cup\A_C$, with $C\subseteq\omega$.

\begin{proposition}
\label{prop:almost-normal-criterion}
For an almost disjoint family $\A$, the following are equivalent.
\begin{enumerate}[label=(\arabic*)]
\item $\Psi(\A)$ is almost-normal.
\item For every $C\subseteq\omega$, there is a partitioner for $\A$ which
separates $\A_C$ from $\A\setminus\A_C$.
\end{enumerate}
\end{proposition}

We shall use the following cardinal characteristics of the continuum.
An \emph{unbounded family} is a set $\mathcal B\subseteq\Baire$ such that for every $f\in\Baire$, there exists $g\in\mathcal B$ with $f\leq^{\infty} g$, i.e., $f(n)\leq g(n)$ for infinitely many $n\in\omega$.
The bounding number is
\[
 \mathfrak b=\min\bigl\{|\mathcal B|:\mathcal B\subseteq\Baire\text{ is unbounded}\bigr\}.
\]

A set $S\subseteq\omega$ \emph{splits} $A\in[\omega]^\omega$ if both
$A\cap S$ and $A\setminus S$ are infinite.
A \emph{splitting family} is a set $\mathcal S\subseteq\Pow(\omega)$ such that every $A\in[\omega]^\omega$ is split by some $S\in\mathcal S$.  The splitting number is
\[
 \mathfrak s=\min\bigl\{|\mathcal S|: \mathcal S\subseteq\Pow(\omega)\text{ is a splitting family}\bigr\}.
\]

\begin{definition}
An \emph{interval partition} of $\omega$ is a sequence
$(I_n)_{n\in\omega}$ of the form $I_n=[t_n,t_{n+1})$, where
$0=t_0<t_1<\cdots$ are natural numbers.
A set $B\subseteq\omega$ \emph{block-splits} this partition if
\[
 |\{n\in\omega:I_n\subseteq B\}|=|\{n\in\omega:I_n\cap B=\emptyset\}|=\omega.
\]
A \emph{block-splitting family} is a set $\mathcal B\subseteq\Pow(\omega)$ such that every interval partition is block-split by some $B\in\mathcal B$.
The \emph{block-splitting number} is
\[
 \mathfrak{s}_{block}=\min\{|\mathcal B|:\mathcal B
 \text{ is a block-splitting family}\}.
\]
\end{definition}
The following result is due to Kamburelis and W\k{e}glorz.
\begin{proposition}[{\cite[Proposition~2.1]{KamburelisWeglorz}}]
\label{prop:block-splitting}
$\mathfrak{s}_{block}=\max\{\mathfrak b,\mathfrak s\}$.
\end{proposition}

A collection $\mathcal C\subseteq[\omega]^\omega$ is \emph{centered} if every finite nonempty subcollection of $\mathcal C$ has an infinite intersection.
A pseudointersection of a centered collection $\mathcal C$ is an infinite set $P\subseteq\omega$ such that $P\subseteq^*C$ for every $C\in\mathcal C$.
The \emph{pseudointersection number} $\mathfrak p$ is the least size of a centered collection
$\F\subseteq[\omega]^\omega$ with no pseudointersection.

Bell's theorem \cite{Bell} states that $\mathfrak p$ is the least cardinal $\kappa$ such that Martin's axiom for $\kappa$ many dense subsets of $\sigma$-centered posets fails; that is, such that there exists a $\sigma$-centered poset $\mathbb P$ and a collection $\mathcal D$ of $\kappa$ many dense subsets of $\mathbb P$ with no filter $G\subseteq\mathbb P$ meeting every $D\in\mathcal D$.
For a proof in this language, see \cite[Theorem~III.3.61]{Kunen}.
It is also well-known that $\mathfrak p\leq\mathfrak b$ (see \cite{Blass}).

For disjoint subfamilies $\B,\C$ of an almost disjoint family, we say
that $\B$ is \emph{weakly separated from $\C$} if some $X\subseteq\omega$
satisfies
\[
 |b\cap X|=\omega\quad(\forall b\in\B)
 \qquad \text{and}\qquad
 |c\cap X|<\omega\quad(\forall c\in\C).
\]
Following Brendle \cite[p.~14]{BrendleDow}, let $\mathfrak{ap}$ be the
least cardinality of an almost disjoint family $\A$ for which some
$\B\subseteq\A$ is not weakly separated from $\A\setminus\B$.

For $t\in 2^{<\omega}$, the basic clopen set generated by $t$ is denoted by $[t]=\{x\in\Cantor:t\subseteq x\}$.

We also use the following generalization of $\lambda$-sets to arbitrary
topological spaces.

\begin{definition}
A topological space $X$ is a \emph{$\lambda$-space} if every countable
subset of $X$ is $G_\delta$ in $X$.  A separable metrizable $\lambda$-space
is called a \emph{$\lambda$-set}.
\end{definition}

Notice that weak $\lambda$-spaces are $T_1$.  For distinct points
$x,y$ of a weak $\lambda$-space $X$, choose a $G_\delta$ set $H$ with
$x\in H$ and $y\notin H$, and write $H=\bigcap_{n\in\omega}U_n$
with each $U_n$ open in $X$.  Some $U_n$ contains $x$ and omits $y$.

We allow $\lambda$-sets and weak $\lambda$-sets to be countable.  This convention
differs from the one in \cite[Definition~1.1]{RodriguesRonchimSzeptycki}, which requires
$\lambda$-sets to be uncountable.
This does not affect the question of whether there exists a weak $\lambda$-set which is not a $\lambda$-set: since every countable set is trivially a $\lambda$-set, any counterexample must be uncountable.

We will use the following notation.

\begin{definition}
    Let
    \begin{itemize}
        \item $\non(\wlambda)=
            \min\bigl\{|X|:X
            \text{ is separable, metrizable, and not a weak $\lambda$-set}\bigr\}$.
        \item $\non(\lambda)=
            \min\bigl\{|X|:X
            \text{ is separable, metrizable, and not a $\lambda$-set}\bigr\}$.
    \end{itemize}
\end{definition}

The following result is due to Rothberger.

\begin{proposition}[{\cite[Lemma~3 and Theorem~1]{Rothberger}}]
\label{prop:non-lambda}
    $\non(\lambda)=\mathfrak b$.
\end{proposition}

For a proof of the preceding result, the reader may also refer to \cite[Theorem~9.1]{vanDouwen}, where $\non(\lambda)$ is denoted by $\mathfrak{b}_\lambda$.

We identify $\Pow(\omega)$ with $\Cantor$ by characteristic functions.
For $x\in\Inf$, let $e_x\in\Baire$ be its increasing enumeration, so
$x=\{e_x(n):n\in\omega\}$.
The map $x\mapsto e_x$ is continuous, since for every $n,m\in\omega$,
the condition $e_x(n)=m$ is determined by $x\cap(m+1)$.

The following criterion is known and is a variation of, e.g., \cite[Lemma~9.3]{vanDouwen} and \cite[Theorem~2.1]{MillerLambdaSets}.
As we did not find a reference for this exact formulation and the proof is short, we include it here for the convenience of the reader.

\begin{lemma}
\label{lem:unbounded-non-lambda}
Let $\Fin\subseteq X\subseteq\Cantor$.
If $\{e_x:x\in X\setminus\Fin\}$ is unbounded in $(\Baire,\leq^*)$,
then $\Fin$ is not $G_\delta$ in $X$.
In particular, $X$ is not a $\lambda$-set.
\end{lemma}

\begin{proof}
Suppose that $\Fin=X\cap G$ for some $G_\delta$ set $G\subseteq\Cantor$.
Then $C=\Cantor\setminus G$ is $F_\sigma$ in the compact space
$\Cantor$, hence $\sigma$-compact, and $X\setminus\Fin\subseteq C\subseteq\Inf$.

Write $C=\bigcup_{j\in\omega}K_j$, where each $K_j$ is compact.
For every $j,n\in\omega$, the set $\{e_x(n):x\in K_j\}$ is finite,
as the continuous image of a compact space in the discrete space $\omega$.
Define $g_j,g\in\Baire$ by:
\[
 g_j(n)=\max\bigl(\{e_x(n):x\in K_j\}\cup\{0\}\bigr)
 \qquad \text{and}\qquad
 g(n)=\max_{j\leq n}g_j(n).
\]
If $x\in K_j$, then $e_x(n)\leq g_j(n)\leq g(n)$ for every $n\geq j$,
so $e_x\leq^*g$.
Thus $\{e_x:x\in X\setminus\Fin\}$ is bounded, contradicting the hypothesis.
Since $\Fin$ is countable, $X$ is not a $\lambda$-set.
\end{proof}

\section{The least cardinality of a non-weak
\texorpdfstring{$\lambda$}{lambda}-set}
\label{sec:weak-lambda}

The question of whether every weak $\lambda$-set is a $\lambda$-set was raised in \cite{RodriguesRonchimSzeptycki}.
In this section we prove that
$\non(\wlambda)=\max\{\mathfrak b,\mathfrak s\}$.
Together with $\non(\lambda)=\mathfrak b$, this gives a weak
$\lambda$-set which is not a $\lambda$-set whenever $\mathfrak b<\mathfrak s$.

We first prove the lower bound using the bounding and splitting numbers.

The following lemma is a well-known characterization of the splitting number $\mathfrak s$.
For a proof, see also \cite[Theorem~3.2]{Blass}.

\begin{lemma}
\label{lem:s-convergence}
The cardinal $\mathfrak s$ is the minimum cardinality of a family
$\mathcal S\subseteq 2^\omega$ such that there is no infinite $A\subseteq\omega$ for which
$r|_A$ converges for every $r\in\mathcal S$.
\end{lemma}

We start proving our characterization.
\begin{lemma}
\label{lem:s-lower}
Every zero-dimensional $T_1$ space $X$ of cardinality less than
$\mathfrak s$ is a weak $\lambda$-space.
\end{lemma}

\begin{proof}
Let $Y,Z\subseteq X$ be disjoint and countable.  If $Y\cup Z$ is finite, then $Y$ is separated from $Z$ by a clopen set.

So assume that $Y\cup Z$ is infinite and enumerate it without repetitions
as $\{d_n:n\in\omega\}$.
Since $X$ is zero-dimensional and $T_1$, for each $n\in\omega$ choose $(U_{n,i}:i\leq n)$ such that:
\begin{enumerate}[label=(\roman*)]
    \item $U_{n,i}$ is a clopen neighborhood of $d_i$ for every $i\leq n$;
    \item For every two distinct $i,j\leq n$, the sets $U_{n,i}$ and $U_{n,j}$ are disjoint.
\end{enumerate}
This is possible as $T_1$ zero-dimensional spaces are Hausdorff.

For each $n\in\omega$, define $\displaystyle C_n=\bigcup\bigl\{U_{n,i}:
 i\leq n\text{ and }d_i\in Y\bigr\}$.
Each $C_n$ is clopen.
For each $x \in X$, let $f_x:\omega\to 2$ be defined by
\[
    f_x(n)=\chi_{C_n}(x)=
    \begin{cases}
    1,&\text{if }x\in C_n,\\
    0,&\text{if }x\notin C_n.
    \end{cases}
\]

In particular, whenever $i\leq n<\omega$:
\begin{equation}
\label{eq:interpolation}
 f_{d_i}(n)=
 \begin{cases}
 1,&\text{if }d_i\in Y,\\
 0,&\text{if }d_i\in Z.
 \end{cases}
\end{equation}

Apply Lemma~\ref{lem:s-convergence} to $\mathcal S=\{f_x:x\in X\}$,
which is possible as $|X|<\mathfrak s$,
to obtain $A\subseteq\omega$ infinite such that $f_x|_A$ converges for every $x\in X$.
Thus, $f_x|A$ is eventually constant for every $x\in X$.
Now let $H$ be the set of all $x\in X$ such that $f_x|_A$ is eventually equal to $1$, that is,
\[
 H=\bigcup_{m\in\omega}\bigcap_{n \in A\setminus m}\{x \in X: f_x(n)=1\}=\bigcup_{m\in\omega}\bigcap_{n \in A\setminus m}C_{n}.
\]
The set $H$ is $F_\sigma$ in $X$. Its complement is
\[
 X\setminus H=\bigcup_{m\in\omega}\bigcap_{n \in A\setminus m}\{x \in X: f_x(n)=0\}=\bigcup_{m\in\omega}\bigcap_{n \in A\setminus m}X\setminus C_{n}.
\]
This is also an $F_\sigma$ set in $X$.
Hence $H$ is also $G_\delta$ in $X$.
Equation~\eqref{eq:interpolation} shows that
$Y\subseteq H$ and
$Z\subseteq X\setminus H$.
\end{proof}

Every metric space of cardinality $<\mathfrak c$ is zero-dimensional.
Indeed, if $(X, d)$ is a metric space with $|X|<\mathfrak c$, then for each $x\in X$, $R_x=\{d(x, y): y \in X\}$ has cardinality $<\mathfrak c$, so $\{B(x, r): r \in (0,\infty)\setminus R_x\}$ is a local base of clopen neighborhoods of $x$.
Thus, we have the following corollary.
\begin{corollary}
    \label{cor:metrizable-weak-lambda}
    Every metrizable space of cardinality less than $\mathfrak s$ is a weak $\lambda$-space.
\end{corollary}

Moreover, we have proved the desired lower bound for $\non(\wlambda)$.

\begin{corollary}
    \label{cor:non-weak-lambda-lower-bound}
    $\max\{\mathfrak b,\mathfrak s\}\leq\non(\wlambda)$.
\end{corollary}

\begin{proof}
As every $\lambda$-set is a weak $\lambda$-set, by Proposition~\ref{prop:non-lambda} we have $\mathfrak b=\non(\lambda)\leq\non(\wlambda)$.
Combining this with Corollary~\ref{cor:metrizable-weak-lambda} gives the desired inequality.
\end{proof}

We now construct a non-weak $\lambda$-set of size
$\max\{\mathfrak b,\mathfrak s\}$ using
Proposition~\ref{prop:block-splitting}.

\begin{theoremproof}{thm:intro-non-weak-lambda}
By Corollary~\ref{cor:non-weak-lambda-lower-bound}, it suffices to prove
$\non(\wlambda)\leq\max\{\mathfrak b,\mathfrak s\}$.
For $\varepsilon\in\{0,1\}$, let:
\[
 Q_\varepsilon=\{x\in\Cantor:x(n)=\varepsilon
 \text{ for all sufficiently large }n\}=\bigcup_{N\in \omega}\bigcap_{n\geq N}\{x\in\Cantor:x(n)=\varepsilon\}.
\]

The sets $Q_0$ and $Q_1$ are disjoint and countable.
By Proposition~\ref{prop:block-splitting}, choose a block-splitting family
$\mathcal B$ of size $\max\{\mathfrak b,\mathfrak s\}$, and put
$ X=Q_0\cup Q_1\cup\{\chi_B:B\in\mathcal B\}\subseteq\Cantor.$

Clearly $|X|=\max\{\mathfrak b,\mathfrak s\}$.
We show that $Q_0,Q_1$ admit no relative $G_\delta$--$F_\sigma$ separator in $X$.

Suppose, towards a contradiction, that a relative
$G_\delta$--$F_\sigma$ set $H\subseteq X$ satisfies
$Q_0\subseteq H$ and $Q_1\subseteq X\setminus H$.
Since both $H$ and $X\setminus H$ are $F_\sigma$ in $X$, write:
\[
 H=X\cap\bigcup_{i\in\omega}F_i,\qquad
 X\setminus H=X\cap\bigcup_{i\in\omega}E_i,
\]
where $F_i,E_i$ are closed in $\Cantor$.
Since $Q_0\subseteq H$ and $Q_1\subseteq X\setminus H$, we have
\[
 F_i\cap Q_1=\varnothing,\qquad E_i\cap Q_0=\varnothing
 \quad(i\in\omega).
\]

Recursively choose $0=t_0<t_1<\cdots$ so that, putting
$I_n=[t_n,t_{n+1})$, for every $x\in\Cantor$,
\begin{equation}
\label{eq:weak-lambda-block}
 \begin{aligned}
 x|_{I_n}\equiv1&\quad\Longrightarrow\quad x\notin \bigcup_{i\leq n}F_i,\\
 x|_{I_n}\equiv0&\quad\Longrightarrow\quad x\notin \bigcup_{i\leq n}E_i.
 \end{aligned}
\end{equation}
To see that this is possible, assume we have defined $t_n$.
We show how to define $t_{n+1}$.

For each $u\in2^{t_n}$,
let $q_u^0=u^\frown\mathbf 0$ and $q_u^1=u^\frown\mathbf 1$,
where $\mathbf 0$ and $\mathbf 1$ are the constant (infinite) sequences with values 0 and 1, respectively.
Since $q_u^1\in Q_1$, we have $q_u^1\notin\bigcup_{i\leq n}F_i$.
Similarly, $q_u^0\notin\bigcup_{i\leq n}E_i$. 
As these unions are closed, choose integers $k_u^0,k_u^1>t_n$ such that
\[
 [q_u^1|_{k_u^1}]\cap\bigcup_{i\leq n}F_i=\varnothing,\qquad
 [q_u^0|_{k_u^0}]\cap\bigcup_{i\leq n}E_i=\varnothing.
\]
There are only finitely many $u\in2^{t_n}$, so we may define
\[
 t_{n+1}=\max\{k_u^\varepsilon:u\in2^{t_n},\ \varepsilon\in\{0,1\}\}.
\]
If $x|_{I_n}\equiv1$, let $u=x|_{t_n}$.
Then $x|_{t_{n+1}}\subseteq q_u^1$.  Since
$k_u^1\leq t_{n+1}$, we have $x\in[q_u^1|_{k_u^1}]$, hence
$x\notin\bigcup_{i\leq n}F_i$.  Similarly, if $x|_{I_n}\equiv0$, then 
$x\notin\bigcup_{i\leq n}E_i$, proving~\eqref{eq:weak-lambda-block}.

Choose $B\in\mathcal B$ that block-splits $(I_n)_{n\in\omega}$,
and put $x=\chi_B$.
Then $x|_{I_n}\equiv0$ for infinitely many $n$ and $x|_{I_n}\equiv1$ for infinitely many $n$.
Fix $i\in\omega$.  Since the set of indices $n$ for which
$x|_{I_n}\equiv1$ is infinite, choose such an $n\geq i$.
By~\eqref{eq:weak-lambda-block},
$x\notin\bigcup_{j\leq n}F_j$, and in particular $x\notin F_i$.
Likewise, $x\notin E_i$.  As $i$ was arbitrary, $x\notin \bigcup_{i\in\omega}(F_i\cup E_i)\supseteq X$, a contradiction.
\end{theoremproof}

\begin{corollary}
\label{cor:b-less-s}
If $\mathfrak b<\mathfrak s$, then there is a weak $\lambda$-set which
is not a $\lambda$-set.
\end{corollary}

\begin{proof}
By Proposition~\ref{prop:non-lambda}, choose a separable metrizable space
$X$ of cardinality $\mathfrak b$ which is not a $\lambda$-set.
Since $|X|=\mathfrak b<\mathfrak s$,
Corollary~\ref{cor:metrizable-weak-lambda} shows that $X$ is a weak
$\lambda$-set.
\end{proof}

\section{A weak \texorpdfstring{$\lambda$}{lambda}-set which is not a
\texorpdfstring{$\lambda$}{lambda}-set under
\texorpdfstring{$\mathfrak p=\mathfrak c$}{p=c}}
\label{sec:direct-ch-weak-lambda}

In this section, we show that $\mathfrak p=\mathfrak c$ implies the
existence of weak $\lambda$-sets that are not $\lambda$-sets. This
complements Corollary~\ref{cor:b-less-s}, which establishes the same
conclusion under $\mathfrak b<\mathfrak s$.

This result was recently obtained by Osipov and Lipin
\cite[Theorem~1.1]{osipovweaklambda}, whose proof uses both Hindman's
finite-unions theorem and Bell's theorem. Our proof, obtained
independently, also uses Hindman's theorem but does not require
Bell's theorem.

In our argument, the finite-unions theorem yields the compact covers
needed for the separators. We first construct a decreasing sequence
of block sequences and then choose points along this sequence.

Throughout this section, we identify $\Pow(\omega)$ with $\Cantor$
via characteristic functions and equip $\Fin$ and $\Inf$ with the
corresponding subspace topologies.

We prove Theorem~\ref{thm:intro-p-c-weak-lambda} by constructing a
witness of cardinality $\mathfrak c$ in $\Cantor$.

\begin{definition}
\label{def:block-sequence}
Let $I\in[\omega]^\omega$.
A \emph{block sequence} indexed by $I$ is a sequence $A=(a_i)_{i\in I}$ of
nonempty finite subsets of $\omega$ such that
$\max a_i<\min a_j$ whenever $i,j\in I$ and $i<j$.
We define:
\begin{enumerate}[label=(\roman*)]
\item the set of all finite unions of blocks of $A$, $\FU(A)=
\left\{\bigcup_{i\in F}a_i:
F\in[I]^{<\omega}\setminus\{\varnothing\}\right\}$;
\item the set of all infinite unions of blocks of $A$,
$[A]_\infty=
\left\{\bigcup_{i\in J}a_i:J\in[I]^\omega\right\}$;
\item the saturation of $A$,
$S(A)=\{r\cup x:r\in[\omega]^{<\omega},\ x\in[A]_\infty\}$.

\end{enumerate}
For $J\in[I]^\omega$, we write $A\restr J=(a_j)_{j\in J}$.  For $m\in\omega$,
we abbreviate $A_{\geq m}=A\restr(I\setminus m)$.
\end{definition}

\begin{definition}
Let $A$ and $B=(b_j)_{j\in J}$ be block sequences.  We say that
$B$ is a \emph{condensation} of $A$, and write $B\preceq A$, if $b_j\in\FU(A)$ for every $j\in J$.

Write $B\preceq^*A$ if $B_{\geq m}\preceq A$ for some $m\in\omega$.
Equivalently, $b_j\in\FU(A)$ for all but finitely many $j\in J$.
\end{definition}

Both $\preceq$ and $\preceq^*$ are easily seen to be transitive.

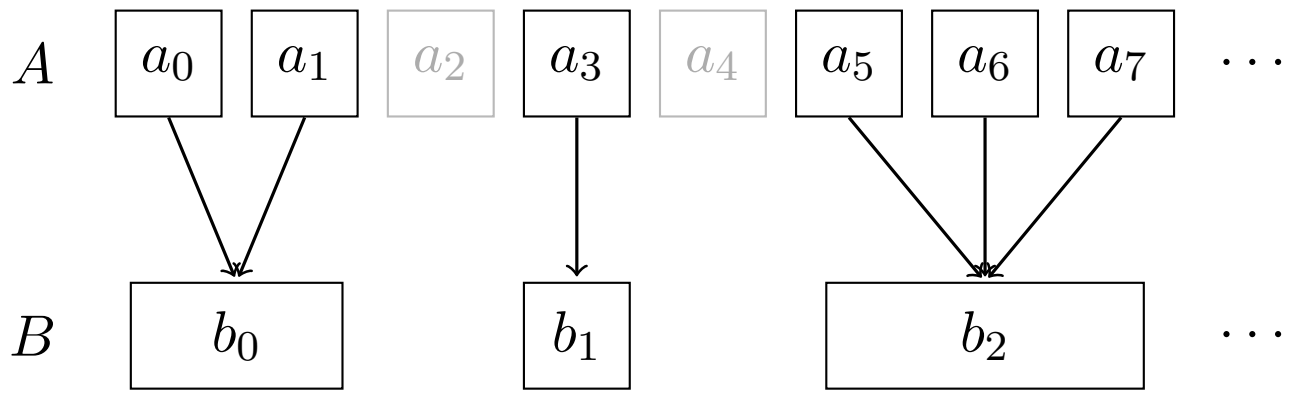
\begin{figure}[H]
\centering
\begin{tikzpicture}[
  block/.style={draw, minimum size=7mm, inner sep=0pt},
  condensed block/.style={draw, minimum height=7mm, inner sep=0pt},
  gap/.style={block, draw=gray!55, text=gray!65},
  block arrow/.style={->, semithick, shorten >=1pt}
]
  \node at (-0.9,0) {$A$};
  \node[block] (a0) at (0,0) {$a_0$};
  \node[block] (a1) at (0.9,0) {$a_1$};
  \node[gap]   (a2) at (1.8,0) {$a_2$};
  \node[block] (a3) at (2.7,0) {$a_3$};
  \node[gap]   (a4) at (3.6,0) {$a_4$};
  \node[block] (a5) at (4.5,0) {$a_5$};
  \node[block] (a6) at (5.4,0) {$a_6$};
  \node[block] (a7) at (6.3,0) {$a_7$};
  \node at (7.2,0) {$\cdots$};

  \node at (-0.9,-1.8) {$B$};
  \node[condensed block, minimum width=14mm] (b0) at (0.45,-1.8) {$b_0$};
  \node[condensed block, minimum width=7mm]  (b1) at (2.7,-1.8) {$b_1$};
  \node[condensed block, minimum width=21mm] (b2) at (5.4,-1.8) {$b_2$};
  \node at (7.2,-1.8) {$\cdots$};

  \draw[block arrow] (a0.south) -- (b0.north);
  \draw[block arrow] (a1.south) -- (b0.north);
  \draw[block arrow] (a3.south) -- (b1.north);
  \draw[block arrow] (a5.south) -- (b2.north);
  \draw[block arrow] (a6.south) -- (b2.north);
  \draw[block arrow] (a7.south) -- (b2.north);
\end{tikzpicture}
\caption{A condensation $B\preceq A$.}
\label{fig:block-condensation}
\end{figure}

Discarding finitely many initial blocks does not change the saturation:
\[
 S(A_{\geq m})=S(A)\qquad(\forall m\in\omega).
\]
Indeed, those blocks can be absorbed into the finite set $r$ in the
definition of $S(A)$.

\begin{lemma}\label{lem:block-lower-bound1}
  Let $A$ and $B$ be block sequences.
  \begin{enumerate}[label=(\roman*)]
    \item If $B\preceq A$, then $\FU(B)\subseteq\FU(A)$.
    \item If $B\preceq A$, then $[B]_\infty\subseteq[A]_\infty$.
    \item If $B\preceq^* A$, then $S(B)\subseteq S(A)$.
  \end{enumerate}
\end{lemma}
\begin{proof}
Write $A=(a_i)_{i\in I}$ and $B=(b_j)_{j\in J}$.

(i) The set $\FU(A)$ is closed under nonempty finite unions.
As each $b_j$ is in $\FU(A)$, every nonempty finite union of blocks of $B$ is in
$\FU(A)$, therefore $\FU(B)\subseteq\FU(A)$.

(ii) Let $x\in[B]_\infty$.
Write $x=\bigcup_{j\in K}b_j$ for some $K\in[J]^\omega$.
For each $j\in K$, choose a nonempty finite $F_j\subseteq I$ such that $b_j=\bigcup_{i\in F_j}a_i$.
Let $H=\bigcup_{j\in K}F_j$.  The sets $F_j$ are pairwise disjoint, so $H$ is infinite.
Then:
\[
  x=\bigcup_{j\in K}b_j
  =\bigcup_{j\in K}\bigcup_{i\in F_j}a_i
  =\bigcup_{i\in H}a_i\in [A]_\infty.
\]

For (iii), choose $m\in\omega$ with $B_{\geq m}\preceq A$.
By (ii), $[B_{\geq m}]_\infty\subseteq[A]_\infty$, so $S(B_{\geq m})\subseteq S(A)$.
Thus it suffices to see that $S(B)\subseteq S(B_{\geq m})$.
Let $x\in S(B)$.
Fix $r\in[\omega]^{<\omega}$ and $K\in[J]^\omega$ with $x=r\cup\bigcup_{j\in K}b_j$.
Then:
\[
  x=r\cup\bigcup_{j\in K}b_j
  =\left(r\cup\bigcup_{j\in K\cap m}b_j\right)\cup\bigcup_{j\in K\setminus m}b_j
  \in S(B_{\geq m}).
\]
\end{proof}

\begin{lemma}\label{lem:block-lower-bound2}
Let $\delta<\mathfrak p$ be an ordinal and let $(A_\xi)_{\xi<\delta}$ be a
$\preceq^*$-decreasing family of block sequences.
Then there exists a block sequence $B$ such that
$B\preceq^*A_\xi$ for every $\xi<\delta$.
\end{lemma}

\begin{proof}
For $m\in\omega$, put $T_m=\{s\in\Fin\setminus\{\varnothing\}:\min s\geq m\}$.
Let
\[
 \mathcal F=\{\FU(A_\xi):\xi<\delta\}
 \cup\{T_m:m\in\omega\}.
\]

Write $N=[\omega]^{<\omega}\setminus\{\varnothing\}$.
Then $\mathcal F$ is a collection of subsets of the infinite countable set $N$.
\begin{claim}
  $\mathcal F$ is centered.
\end{claim}
\begin{proof}[Proof of the claim]
The sets $T_m$ are decreasing for $m\in\omega$ and clearly infinite.
Thus it suffices to show that given a nonempty finite $H\subseteq\delta$ and $m\in\omega$,
\[
  T_m\cap\bigcap_{\xi\in H}\FU(A_\xi)\text{ is infinite.}
\]
Let $\eta=\max H$.
For every $\xi\in H$, choose $m_\xi$ with $(A_\eta)_{\geq m_\xi}\preceq A_\xi$.
Choose $M\geq\max\{m_\xi:\xi\in H\}$ such that every block of
$(A_\eta)_{\geq M}$ has minimum at least $m$.
Then $(A_\eta)_{\geq M}\preceq A_\xi$ for every $\xi\in H$.

Thus, $\FU((A_\eta)_{\geq M})\subseteq\FU(A_\xi)$ for every $\xi\in H$.
Also, given $x\in\FU((A_\eta)_{\geq M})$, we have $\min x\geq m$, so $x\in T_m$.
As $\FU((A_\eta)_{\geq M})$ is infinite, this proves the claim.
\end{proof}

Since $|\mathcal F|\leq|\delta|+\aleph_0<\mathfrak p$, choose an infinite pseudointersection
$P$ of $\mathcal F$.
Since $P\subseteq^*T_m$ for every $m$, we can recursively choose
$b_n\in P$ with
$\max b_n<\min b_{n+1}$.
Put $B=(b_n)_{n\in\omega}$.
For each $\xi<\delta$, the inclusion $P\subseteq^*\FU(A_\xi)$ implies
that all but finitely many $b_n$ belong to $\FU(A_\xi)$.
Thus $B\preceq^*A_\xi$.
\end{proof}

We use Hindman's finite-unions theorem, which we state below.

\begin{proposition}[{\cite[Corollary~3.3]{Hindman}}]
\label{prop:hindman-finite-unions}
For every block sequence $A$ indexed by a set $I\in[\omega]^\omega$ and every finite coloring
$d:\FU(A)\to L$, where $L$ is a finite set, there exists a block sequence
$B\preceq A$ indexed by $I$ such that $d$ is constant on $\FU(B)$.
\end{proposition}

For the next lemma, we need an \emph{almost disjoint family} of cardinality $\mathfrak c$.
An almost disjoint family on an infinite countable set $N$ is a collection $\mathcal A\subseteq[N]^\omega$ such that $|a\cap b|<\omega$ for every distinct $a,b\in\mathcal A$.
It is well-known that there are almost disjoint families of cardinality $\mathfrak c$.
For instance, for each $z\in\Cantor$, let $a_z=\{z|n: n \in \omega\}$.
Then $\{a_z: z \in \Cantor\}$ is an almost disjoint family of cardinality $\mathfrak c$ on $2^{<\omega}$.
A bijection of $2^{<\omega}$ with $N$ gives an almost disjoint family of cardinality $\mathfrak c$ on $N$.

\begin{lemma}
\label{lem:avoid-small-set}
Let $A$ be a block sequence and let $E\subseteq\Inf$
have cardinality less than $\mathfrak c$.  There exists a subsequence $B$ of $A$ such that
$S(B)\cap E=\varnothing$.
\end{lemma}

\begin{proof}
Write $A=(a_i)_{i\in I}$, where $I\in[\omega]^\omega$.
Let $\mathcal A\subseteq[I]^\omega$ be an almost disjoint family of cardinality $\mathfrak c$.
For distinct $a,b\in\mathcal A$, the unions $\bigcup_{i\in a}a_i$ and
$\bigcup_{i\in b}a_i$ have finite intersection, so
$S(A\restr a)\cap S(A\restr b)=\varnothing$.
Consequently, $|\{a\in\mathcal A:S(A\restr a)\cap E\neq\varnothing\}|\leq|E|<\mathfrak c$.
Choose $a\in\mathcal A$ with $S(A\restr a)\cap E=\varnothing$ and put $B=A\restr a$.
\end{proof}

The next lemma says that whenever $[\omega]^{<\omega}$ is colored with two colors,
every block sequence has a condensation whose saturation is covered by countably many compact sets whose intersections with $[\omega]^{<\omega}$ are nonempty and monochromatic.

\begin{lemma}
\label{lem:block-fusion}
Let $A$ be a block sequence and let $c:\Fin\to2$.
There exist $B\preceq A$ and a countable family $\mathcal K$ of compact
subsets of $\mathcal P(\omega)\approx\Cantor$ such that:
\begin{enumerate}[label=(\Roman*)]
\item $S(B)\subseteq\bigcup\mathcal K$;
\item for every $K\in\mathcal K$, $K\cap\Fin\neq\varnothing$;
\item for every $K\in\mathcal K$, $c|_{K\cap\Fin}$ is constant;
\item for every $K\in\mathcal K$, $K\cap\Inf\subseteq S(B)$.
\end{enumerate}
\end{lemma}

\begin{proof}
Enumerate $\Fin=\{z_j:j\in\omega\}$.
We recursively construct sequences $(R_n)_{n\in\omega}$ and $(b_n)_{n\in\omega}$ such that, for every
$n\in\omega$:
\begin{enumerate}[label=(\roman*)]
  \item $R_n$ is a block sequence;
  \item $b_n\in\FU(R_n)$;
  \item $R_0=A$;
  \item $\max b_n<\min b_{n+1}$;
  \item $R_{n+1}\preceq R_n$, and $\max b_n<\min s$
  for every block $s$ of $R_{n+1}$;
  \item For every $j\leq n$ and $F\subseteq n$, the map
  \[
    v\longmapsto
    c\left(z_j\cup\bigcup_{i\in F}b_i\cup b_n\cup v\right)
  \]
  is constant on $\FU(R_{n+1})\cup\{\varnothing\}$.
\end{enumerate}
Start with $R_0=A$.
Assume that $(R_i:i\leq n)$ and $(b_i:i<n)$ have been chosen.
We show how to choose $b_n$ and $R_{n+1}$.
Put $Q_n=(n+1)\times\mathcal P(n)$ and define
\[
 d_n:\FU(R_n)\longrightarrow 2^{Q_n}\qquad \text{ and } \qquad 
 d_n(v)=
 \left(c\left(z_j\cup\bigcup_{i\in F}b_i\cup v\right)\right)
 _{(j,F)\in Q_n}.
\]
Let $I$ be the index set of $R_n$.
Since $2^{Q_n}$ is finite,
apply Hindman's finite-unions theorem (Proposition~\ref{prop:hindman-finite-unions}) to $R_n$ and $d_n$
to obtain a condensation $H_n=(h_i)_{i\in I}\preceq R_n$ such that
$d_n$ is constant on $\FU(H_n)$.
Put $i_*=\min I$ and define
\[
  b_n=h_{i_*}\qquad \text{ and } \qquad 
  R_{n+1}=H_n\restr\bigl(I\setminus\{i_*\}\bigr).
\]
Clearly, all the required properties are satisfied.

Now let $B=(b_n)_{n\in\omega}$. Then $B$ is a block sequence, $B\preceq A$, and
$B_{\geq n+1}\preceq R_{n+1}$ for every $n\in \omega$.
For $j\leq n$ and $F\subseteq n$, put
\begin{equation}
\label{eq:fusion-compact}
 K_{j,F,n}=\left\{
 z_j\cup\bigcup_{i\in F}b_i\cup b_n\cup\bigcup_{k\in J}b_k:
 J\subseteq\{k\in\omega:k>n\}\right\}.
\end{equation}
\begin{claim}
  Each set $K_{j,F,n}$ is compact.
\end{claim}
\begin{proof}[Proof of the claim]
Fix $n\in\omega$, $j\leq n$, and $F\subseteq n$.  Put
$T_n=\{k\in\omega:k>n\}$ and
$u=z_j\cup\bigcup_{i\in F}b_i\cup b_n$.
Define
\[
 \Phi_{j,F,n}:2^{T_n}\longrightarrow\mathcal P(\omega)\approx\Cantor,\qquad
 \Phi_{j,F,n}(t)=u\cup\bigcup_{\substack{k\in T_n\\t(k)=1}}b_k.
\]
The range of $\Phi_{j,F,n}$ is $K_{j,F,n}$.
We check that $\Phi_{j,F,n}$ is continuous, which then guarantees that $K_{j,F,n}$ is compact.
For brevity, write $\Phi=\Phi_{j,F,n}$.  The sets
\[
 V_m^1=\{x\subseteq\omega:m\in x\}\qquad \text{ and } \qquad 
 V_m^0=\{x\subseteq\omega:m\notin x\},
\]
for $m \in \omega$
form a subbasis for the topology of $\mathcal P(\omega)$.
As for each $m \in \omega$, $V_m^0=\mathcal P(\omega)\setminus V_m^1$, it suffices to check that the preimages of the sets $V_m^1$ are clopen sets.

\textbf{Case 1:} $m\in u$.  Then $\Phi(t)\in V_m^1$ for every $t\in2^{T_n}$, so $\Phi^{-1}[V_m^1]=2^{T_n}$, which is clopen.

\textbf{Case 2:} $m\in b_k\setminus u$ for some $k\in T_n$. We claim that $\Phi^{-1}[V_m^1]=\{t\in2^{T_n}:t(k)=1\}$, which is clopen.
If $t(k)=1$, then $b_k\subseteq\Phi(t)$, so $m\in\Phi(t)$.
Conversely, if $m\in\Phi(t)$, then $m\notin u$ implies that there exists $k'\in T_n$ with $t(k')=1$ and $m\in b_{k'}$.  Since the blocks are pairwise disjoint, $k=k'$, so $t(k)=1$.

\textbf{Case 3:} $m\notin u\cup\bigcup_{k\in T_n}b_k$.  Then $m\notin\Phi(t)$ for every $t\in2^{T_n}$, so $\Phi^{-1}[V_m^1]=\varnothing$, which is clopen.

Thus the preimage of every subbasic open set is open, so $\Phi$ is
continuous.  Since $2^{T_n}$ is compact, its image $K_{j,F,n}$ is compact.
\end{proof}

Put $\mathcal K=\{K_{j,F,n}:n\in\omega,\ j\leq n,\ F\subseteq n\}$.

$\mathcal K$ is countable, and, by the Claim, each element of $\mathcal K$ is compact.
We verify (I)--(IV).

(I) Let $x\in S(B)$ and write $x=r\cup\bigcup_{k\in J}b_k$,
where $r\in\Fin$ and $J\in[\omega]^\omega$.
Choose $j$ with $r=z_j$.  Since $J$ is infinite, choose $n\in J$
with $n\geq j$, and put $F=J\cap n$ and $L=J\setminus(n+1)$.
Then $F\subseteq n$, $L\subseteq\{k\in\omega:k>n\}$, and
\[
 x=z_j\cup\bigcup_{i\in F}b_i\cup b_n\cup\bigcup_{k\in L}b_k
 \in K_{j,F,n}.
\]
Thus $S(B)\subseteq\bigcup\mathcal K$.

For the remaining properties, fix $n\in\omega$, $j\leq n$, and
$F\subseteq n$, and put $u=z_j\cup\bigcup_{i\in F}b_i\cup b_n$.

(II) Taking $J=\varnothing$ in \eqref{eq:fusion-compact} gives
$u\in K_{j,F,n}\cap\Fin$, so this intersection is nonempty.

(III) Let $x\in K_{j,F,n}\cap\Fin$ and write
$x=u\cup\bigcup_{k\in J}b_k$, where $J\subseteq\{k\in\omega:k>n\}$.
Since the blocks are nonempty and pairwise disjoint, $J$ is finite.
The relation $B_{\geq n+1}\preceq R_{n+1}$ implies that
\[
 v=\bigcup_{k\in J}b_k\in\FU(R_{n+1})\cup\{\varnothing\}.
\]
As $x=u\cup v$, we have $c(x)=c(u)$
by property (vi).
Hence $c\restr(K_{j,F,n}\cap\Fin)$ is constant.

(IV) Let $x\in K_{j,F,n}\cap\Inf$ and write
$x=u\cup\bigcup_{k\in J}b_k$, where $J\subseteq\{k\in\omega:k>n\}$.
Since $u$ and all the blocks are finite, $J$ is infinite.  Therefore
\[
 x=z_j\cup\bigcup_{k\in F\cup\{n\}\cup J}b_k\in S(B).
\]
This proves $K_{j,F,n}\cap\Inf\subseteq S(B)$ and completes the proof.

\end{proof}

\begin{lemma}
\label{lem:block-unbounded}
For every block sequence $A$, the set $\{e_x:x\in[A]_\infty\}$
is unbounded in $(\Baire,\leq^*)$.
\end{lemma}

\begin{proof}
Write $A=(a_i)_{i\in I}$ and fix $g\in\Baire$.
Choose $(i_k:k\in\omega)$ and $(m_k:k\in\omega)$ such that, for every $k\in\omega$:
\begin{enumerate}[label=(\roman*)]
  \item $i_k \in I$;
  \item $i_{k+1}>i_k$;
  \item $m_k=\sum_{j<k}|a_{i_j}|$;
  \item $\min a_{i_k}>g(m_k)$.
\end{enumerate}
Such choices are possible because the minima of the blocks tend to infinity.

Put $x=\bigcup_{k\in\omega}a_{i_k}\in[A]_\infty$.
Given $k \in \omega$, by construction, exactly $m_k$ elements of $x$ precede $\min a_{i_k}$,
so $e_x(m_k)=\min a_{i_k}>g(m_k)$.
Thus $g<^{\infty} e_x$, as required.
\end{proof}

We also record the following elementary property of $F_\sigma$ sets.

\begin{lemma}
\label{lem:fsigma-reduction}
Let $M$ be a metrizable space and let $P_0,P_1$ be $F_\sigma$ subsets of
$M$ with $M=P_0\cup P_1$.
There exist disjoint $F_\sigma$ subsets $H_0,H_1$ of $M$ such that
\[
 M=H_0\cup H_1\qquad \text{and} \qquad H_i\subseteq P_i\quad(i\in\{0,1\}).
\]
\end{lemma}

\begin{proof}
For each $i\in\{0,1\}$, write $P_i=\bigcup_{n\in\omega}C_{i,n}$,
where $(C_{i,n})_{n\in\omega}$ is an increasing sequence of closed
subsets of $M$.
We may assume that $C_{0, 0}=C_{1, 0}=\varnothing$.
Define
\[
 H_0=\bigcup_{n\in\omega}(C_{0,n+1}\setminus C_{1,n})\qquad \text{and} \qquad
 H_1=\bigcup_{n\in\omega}(C_{1,n+1}\setminus C_{0,n+1}).
\]
Each difference of closed sets is $F_\sigma$ in the metrizable space
$M$, so both $H_i$ are $F_\sigma$, and clearly $H_i\subseteq P_i$.
It is now straightforward to check that $H_0\cap H_1=\varnothing$ and $M=H_0\cup H_1$.
\end{proof}

Combining Lemma~\ref{lem:fsigma-reduction} with Proposition~\ref{prop:non-lambda} gives the
following criterion.

\begin{lemma}
\label{lem:weak-lambda-cover}
Let $X$ be a separable metrizable space.
Suppose that, for every pair of disjoint countable subsets $Y,Z\subseteq X$,
there is a countable family $\mathcal F$ of closed subsets of $X$ such that
\begin{itemize}
  \item $\left|X\setminus\bigcup\mathcal F\right|<\mathfrak b$;
  \item for every $F\in\mathcal F$, either $F\cap Y=\varnothing$ or $F\cap Z=\varnothing$.
\end{itemize}
Then $X$ is a weak $\lambda$-set.
\end{lemma}

\begin{proof}
Fix disjoint countable subsets $Y,Z\subseteq X$ and choose
$\mathcal F$ as in the hypothesis.  Put
\[
 T=\left(X\setminus\bigcup\mathcal F\right)\cup Y\cup Z.
\]
Since $|T|<\mathfrak b$, Proposition~\ref{prop:non-lambda} implies that
$T$ is a $\lambda$-set.  Thus $Y=T\cap G$ for some $G_\delta$ subset
$G$ of $X$.  Define
\begin{align*}
  P_0&=Y\cup\bigcup\{F\in\mathcal F:F\cap Z=\varnothing\},\\
  P_1&=(X\setminus G)\cup
  \bigcup\{F\in\mathcal F:F\cap Y=\varnothing\}.
\end{align*}
Both sets are $F_\sigma$ in $X$.
Since $T\subseteq Y\cup(X\setminus G)$ and every member of
$\mathcal F$ is contained in $P_0$ or $P_1$, it follows that $\{P_0,P_1\}$ is a cover of $X$ in which $P_0$ and $P_1$ are $F_\sigma$ sets.

Apply Lemma~\ref{lem:fsigma-reduction} to obtain a partition
$X=H_0\cup H_1$ into $F_\sigma$ sets with $H_i\subseteq P_i$.
Then $Y\cap H_1=\varnothing$ and $Z\cap H_0=\varnothing$, so $Y\subseteq H_0$ and $Z\subseteq H_1$.
Each $H_i$ is also $G_\delta$, since its complement is $F_\sigma$.
Hence $X$ is a weak $\lambda$-set.
\end{proof}

We are now ready to prove Theorem~\ref{thm:intro-p-c-weak-lambda}.

\begin{theoremproof}{thm:intro-p-c-weak-lambda}
Assume $\mathfrak p=\mathfrak c$.
Enumerate $2^{\Fin}\times[\Inf]^{\leq\omega}$ as $\{(c_\alpha,E_\alpha):\alpha<\mathfrak c\}$.

First recursively construct sequences $(A_\alpha)_{\alpha<\mathfrak c}$
and $(\mathcal K_\alpha)_{\alpha<\mathfrak c}$ such that, for every
$\alpha<\mathfrak c$:
\begin{enumerate}[label=(\roman*),ref=(\roman*)]
  \item\label{item:ch-block-sequence} $A_\alpha$ is a block sequence;
  \item\label{item:ch-block-tower} $A_\alpha\preceq^*A_\xi$ for every $\xi<\alpha$;
  \item\label{item:ch-compact-family} $\mathcal K_\alpha$ is a countable family of compact subsets of $\Cantor$;
  \item\label{item:ch-avoid} $S(A_\alpha)\cap E_\alpha=\varnothing$;
  \item\label{item:ch-cover} $S(A_\alpha)\subseteq\bigcup\mathcal K_\alpha$;
  \item\label{item:ch-compact-properties} for every $K\in\mathcal K_\alpha$,
  \[
    \varnothing\neq K\cap\Fin,
    \qquad c_\alpha\restr(K\cap\Fin)\text{ is constant},
    \qquad K\cap\Inf\subseteq S(A_\alpha).
  \]
\end{enumerate}

Assume that $(A_\xi)_{\xi<\alpha}$ and $(\mathcal K_\xi)_{\xi<\alpha}$
have been chosen.  We show how to choose $A_\alpha$ and $\mathcal K_\alpha$.
If $\alpha=0$, let $D=(\{n\})_{n\in\omega}$.
If $\alpha>0$, apply
Lemma~\ref{lem:block-lower-bound2} to obtain a block sequence
$D\preceq^*A_\xi$ for every $\xi<\alpha$.
This is possible because
the previously chosen block sequences are $\preceq^*$-decreasing
and $\alpha<\mathfrak p$.

Use Lemma~\ref{lem:avoid-small-set} to choose a subsequence $C$
of $D$ such that $S(C)\cap E_\alpha=\varnothing$.
Apply Lemma~\ref{lem:block-fusion} to $C$ and $c_\alpha$ to obtain
a block sequence $A_\alpha\preceq C$ and a family $\mathcal K_\alpha$
satisfying \ref{item:ch-compact-family}, \ref{item:ch-cover}, and
\ref{item:ch-compact-properties}.
Since $A_\alpha\preceq C\preceq D$, transitivity gives \ref{item:ch-block-tower}.
By Lemma~\ref{lem:block-lower-bound1}, $S(A_\alpha)\subseteq S(C)$,
so \ref{item:ch-avoid} holds as well.
This ends the recursive construction.

This fixes the block sequences and covers independently of the points.
Now enumerate $\Baire=\{g_\alpha:\alpha<\mathfrak c\}$.
For each $\alpha<\mathfrak c$, use Lemma~\ref{lem:block-unbounded} to choose $x_\alpha\in[A_\alpha]_\infty$ such that $g_\alpha<^{\infty}e_{x_\alpha}$.
Put
\[
 X=\Fin\cup\{x_\alpha:\alpha<\mathfrak c\}.
\]

Thus,
$\{e_{x_\alpha}:\alpha<\mathfrak c\}$ is unbounded in
$(\Baire,\leq^*)$.  By the definition of $\mathfrak b$, this collection
has cardinality at least $\mathfrak b=\mathfrak c$, so $|X|=\mathfrak c$.
We now show that $X$ is a weak $\lambda$-set but not a $\lambda$-set.

\smallskip
\noindent\emph{$X$ is a weak $\lambda$-set.}
Let $Y,Z\subseteq X$ be countable and disjoint.
Put $E=(Y\cup Z)\cap\Inf$ and choose a coloring $c:\Fin\to2$ with
value zero on $Y\cap\Fin$ and one on $Z\cap\Fin$.
Choose $\alpha$ with $(c,E)=(c_\alpha,E_\alpha)$.  For every $\beta\geq\alpha$,
Property \ref{item:ch-block-tower} and the choice of $x_\beta$, together
with Lemma~\ref{lem:block-lower-bound1}, give
\begin{equation}
\label{eq:ch-future-in-tail}
 x_\beta\in[A_\beta]_\infty\subseteq S(A_\beta)\subseteq S(A_\alpha).
\end{equation}
Also by \ref{item:ch-compact-properties}, for every $K\in\mathcal K_\alpha$, $c_\alpha\restr(K\cap\Fin)$ is constant.
Let $\varepsilon(K)\in\{0,1\}$ be the
unique value such that $c_\alpha\restr(K\cap\Fin)\equiv\varepsilon(K)$.

By the choice of $c_\alpha=c$, $K$ avoids $Z\cap\Fin$ if
$\varepsilon(K)=0$, and avoids $Y\cap\Fin$ if $\varepsilon(K)=1$.
Moreover, $K\cap\Inf\subseteq S(A_\alpha)$ and
$S(A_\alpha)\cap E_\alpha=\varnothing$, so every $K\in\mathcal K_\alpha$
avoids $(Y\cup Z)\cap\Inf=E_\alpha$.
Thus, for every $K\in\mathcal K_\alpha$,
\begin{equation}
\label{eq:ch-compact-avoidance}
 \varepsilon(K)=0\Longrightarrow K\cap Z=\varnothing
 \qquad \text{and} \qquad
 \varepsilon(K)=1\Longrightarrow K\cap Y=\varnothing.
\end{equation}

Put $\mathcal F=\{K\cap X:K\in\mathcal K_\alpha\}$.
By Property~\ref{item:ch-compact-family}, $\mathcal F$ is countable,
and each of its members is closed in $X$.
Equation~\eqref{eq:ch-compact-avoidance} shows that each member of
$\mathcal F$ is disjoint from $Y$ or from $Z$.
By \eqref{eq:ch-future-in-tail} and Property~\ref{item:ch-cover},
\[
 X\setminus\bigcup\mathcal F
 \subseteq\Fin\cup\{x_\xi:\xi<\alpha\}.
\]
The set on the right has cardinality less than $\mathfrak b=\mathfrak c$.
Since $Y,Z$ were arbitrary, Lemma~\ref{lem:weak-lambda-cover} implies
that $X$ is a weak $\lambda$-set.

\smallskip
\noindent\emph{$X$ is not a $\lambda$-set.}
By Lemma~\ref{lem:unbounded-non-lambda}, the unboundedness of
$\{e_{x_\alpha}:\alpha<\mathfrak c\}=\{e_x:x\in X\setminus\Fin\}$
implies that $\Fin$ is not $G_\delta$ in $X$.
Hence $X$ is not a $\lambda$-set.
\end{theoremproof}

We emphasize that the proof presented in this section is direct
and does not rely on Bell's theorem, in contrast to
the approach adopted in \cite{osipovweaklambda}.

\section{Weak \texorpdfstring{$\lambda$}{lambda}-sets when
\texorpdfstring{$\mathfrak s<\mathfrak b$}{s < b}}
\label{sec:splitting-bounding-obstruction}

The inequality $\mathfrak b<\mathfrak s$ yields a weak $\lambda$-set
which is not a $\lambda$-set.  When $\mathfrak s<\mathfrak b$,
the two classes coincide.  The proof applies to arbitrary first countable
spaces.

Recall that every weak $\lambda$-space is $T_1$.

\begin{proposition}
\label{prop:s-b-first-countable}
Assume $\mathfrak s<\mathfrak b$.
Every first countable weak $\lambda$-space is a $\lambda$-space.
\end{proposition}

\begin{proof}
Let $X$ be a first countable weak $\lambda$-space.
Since $X$ is $T_1$, first countability makes every singleton, and hence every
finite subset, $G_\delta$ in $X$.
Let $D\subseteq X$ be countably infinite.
For each $d\in D$, fix a decreasing open local base
$(U(d,n): n \in \omega)$.  Since $X$ is $T_1$, for every $d\in D$ we have
\begin{equation}
\label{eq:s-b-local-bases}
 \bigcap_{n\in\omega}U(d,n)=\{d\}
\end{equation}

Fix a splitting family of subsets of $D$, 
$\mathcal S\subseteq\Pow(D)$, of cardinality $\mathfrak s$.
For $S\in\mathcal S$, write $S^1=S$ and $S^0=D\setminus S$.
For each $S\in\mathcal S$, choose a partition
$X=H_S^0\mathbin{\dot\cup}H_S^1$ into $F_\sigma$ sets such that $S^i\subseteq H_S^i$ for $i<2$.  This is possible because $X$ is a weak $\lambda$-space.
For $\varepsilon<2$, write
\[
 H_S^\varepsilon=\bigcup_{n\in\omega}F_{S,n}^\varepsilon,
\]
where each $F_{S,n}^\varepsilon$ is closed in $X$ and the sequence is
increasing in $n$.

Choose an increasing sequence $(D_n)_{n\in\omega}$
of finite subsets of $D$ with $D=\bigcup_n D_n$.

For each $S\in\mathcal S$, choose $q_S\in\Baire$ such that,
for every $n\in\omega$,
\begin{equation}
\label{eq:s-b-neighborhoods}
 \begin{aligned}
 U(d,q_S(n))\cap F_{S,n}^0&=\varnothing
 &&\text{for }d\in D_n\cap S^1,\\
 U(d,q_S(n))\cap F_{S,n}^1&=\varnothing
 &&\text{for }d\in D_n\cap S^0.
 \end{aligned}
\end{equation}
Such a choice is possible because $D_n$ is finite, $F_{S,n}^0$ avoids
$S^1$, $F_{S,n}^1$ avoids $S^0$, and the local bases are decreasing.

As $\mathfrak s<\mathfrak b$, $\{q_S:S\in\mathcal S\}$ is bounded in $(\Baire,\leq^*)$.
Choose $h\in\Baire$ such that
$q_S\leq^*h$ for every
$S\in\mathcal S$.
We may assume that for every $n\in\omega$, $h(n)\ge n$.
Define
\begin{equation}
\label{eq:s-b-gdelta}
 G=\bigcap_{m\in\omega}\bigcup_{n\ge m}\bigcup_{d\in D_n}
 U(d,h(n)).
\end{equation}
This is a $G_\delta$ subset of $X$, and $D\subseteq G$.
We will show that $D=G$.

Suppose towards a contradiction that $G\setminus D\neq \varnothing$ and fix $x\in G\setminus D$.
By Equation~\eqref{eq:s-b-gdelta}, there exists an infinite set
$N\subseteq\omega$ and points $d_n\in D_n$, for $n\in N$, such that
\[
 x\in \bigcap_{n \in N}U(d_n,h(n)).
\]
For each $d\in D$, Equation~\eqref{eq:s-b-local-bases} gives some
$M'$ with $x\notin U(d,M')$.
If $n\in N$ and $n\ge M'$, then $h(n)\ge n\ge M'$, so
$U(d,h(n))\subseteq U(d,M')$.
Hence $x\notin U(d,h(n))$, whereas $x\in U(d_n,h(n))$,
and therefore $d_n\ne d$.
Thus each point of $D$ occurs only finitely often among the $d_n$, which implies that $A=\{d_n:n\in N\}\subseteq D$ is infinite.

Choose $S\in\mathcal S$ which splits $A$, and let $\varepsilon<2$ be such that
$x\in H_S^\varepsilon$.
Since the closed sets $F_{S,n}^\varepsilon$ increase to
$H_S^\varepsilon$ and $q_S\leq^*h$, choose $M\in\omega$ such that
$x\in F_{S,M}^\varepsilon$ and $q_S(n)\le h(n)$ for every $n\ge M$.
For $n\in N$ with $n\ge M$, we have $x\in F_{S,n}^\varepsilon$ and
\begin{equation}
\label{eq:s-b-witness-neighborhood}
    x\in U(d_n,h(n))\subseteq U(d_n,q_S(n)).
\end{equation}
Equation~\eqref{eq:s-b-neighborhoods} then gives
$d_n\in S^\varepsilon$: otherwise $U(d_n,q_S(n))$ would be
disjoint from $F_{S,n}^\varepsilon$, contradicting
Equation~\eqref{eq:s-b-witness-neighborhood}.
Consequently, $A\setminus S^\varepsilon\subseteq\{d_n:n\in N,\ n<M\}$,
contradicting that $S$ splits $A$.

Therefore $G=D$.  Since $D$ was arbitrary, $X$ is a $\lambda$-space.
\end{proof}

\begin{theoremproof}{thm:intro-s-b-weak-lambda}
Every weak $\lambda$-set is first countable, so the conclusion follows
from Proposition~\ref{prop:s-b-first-countable}.
\end{theoremproof}

\section{Weak \texorpdfstring{$\lambda$}{lambda}-sets in the Cohen model}
\label{sec:cohen-model}

In the previous sections, we saw that both $\mathfrak p=\mathfrak c$ and
$\mathfrak b<\mathfrak s$ imply the existence of weak $\lambda$-sets
which are not $\lambda$-sets.
Proposition~\ref{prop:s-b-first-countable} gives the opposite conclusion when
$\mathfrak s<\mathfrak b$.  In this section, we obtain the same
conclusion by adding $\aleph_2$ Cohen reals to a model of $\CH$.

For a set $A$, let $\Fn(A,2)$ denote the set of finite partial
functions from $A$ into $2$, ordered by reverse inclusion:
$p\leq q$ if and only if $p\supseteq q$.  For a cardinal $\kappa$,
the Cohen forcing for adding $\kappa$ many Cohen reals is $\mathbb C_\kappa=\Fn(\kappa\times\omega,2)$.
The Cohen forcing for adding a single Cohen real is $\mathbb C=\Fn(\omega,2)$.

To prove Theorem~\ref{thm:intro-cohen-weak-lambda}, we first determine
how adding one Cohen real affects whether a
ground-model set is a $\lambda$-set or a weak $\lambda$-set.
We will need the following well-known result.

\begin{lemma}[{\cite[Corollary~2.2]{Goldstern}}]
\label{cohen:countable-cover}
Let $\mathbb P$ be a ccc forcing, and let $X$ be a set.
Suppose that $p\in\mathbb P$ and $\dot D$ is a $\mathbb P$-name such that
\[
 p\Vdash\dot D\subseteq\check X
             \text{ and }|\dot D|\leq\aleph_0.
\]
Then there exists a countable $A\subseteq X$ such that $p\Vdash\dot D\subseteq\check A$.
\end{lemma}

\begin{lemma}
\label{cohen:coloring}
Let $X\subseteq\Cantor$, and let $D\subseteq X$ be countable.
If $D$ is not $G_\delta$ in $X$, then
\[
 \mathbb C\Vdash
 \begin{gathered}
 \exists\,B\subseteq\check D\,\left(\text{no relatively }
 F_\sigma\text{ and }G_\delta\text{ subset }H\text{ of }\check X
 \text{ satisfies }H\cap\check D=B\right).
 \end{gathered}
\]
\end{lemma}

\begin{proof}
Work over a countable transitive model $M$ of $\ZFC$ with $X,D\in M$
satisfying the hypotheses there.
Let $G'\subseteq\mathbb C$ be generic over $M$.
Since $D$ is countably infinite in $M$, $\Fn(\omega,2)$ and
$\Fn(D,2)$ are isomorphic in $M$.  Thus there exists a
$\Fn(D,2)$-generic filter $G$ over $M$ such that $M[G]=M[G']$.

Let $g=\bigcup G:D\to 2$.
Then $g \in M[G]$.
Let $B_0=g^{-1}[\{0\}]$ and $B_1=g^{-1}[\{1\}]$.

Working in $M[G]$, suppose that there exists $H\subseteq X$ which
is $F_\sigma$ and $G_\delta$ in $X$ and satisfies $H\cap D=B_1$.
Choose $G_\delta$ sets $H_0,H_1\subseteq\Cantor$ such that
\[
 H_0\cap X=X\setminus H,\qquad H_1\cap X=H.
\]
For each $i<2$, choose a sequence $U^i=(U_n^i:n\in\omega)$ of open
subsets of $\Cantor$ such that $H_i=\bigcap_{n\in\omega}U_n^i$.

Let $\dot U^i\in M$ be a $\Fn(D, 2)$-name such that $\dot U^i_G=U^i$ for each $i<2$.
By the Maximum Principle in $M$, let $\dot g$ be a $\Fn(D, 2)$-name such that $\Fn(D, 2)\Vdash\dot g=\bigcup\dot G$.

By the Truth Lemma, there exists $p_0\in G$ such that
\begin{enumerate}[label=(\roman*)]
      \item $p_0\Vdash\dot U^0 \text{ and }\dot U^1\text{ are sequences of open subsets of }\Cantor$;
      \item \label{item:intersection} $p_0\Vdash\check X\cap\bigcap_{n\in\omega}\dot U^0(n)
          \cap\bigcap_{n\in\omega}\dot U^1(n)=\varnothing$;
      \item $p_0\Vdash\forall d\in\check D\,\forall n\in\omega\,
          \bigl(\dot g(d)=0\Rightarrow d\in\dot U^0(n)\bigr)$.
      \item $p_0\Vdash\forall d\in\check D\,\forall n\in\omega\,
          \bigl(\dot g(d)=1\Rightarrow d\in\dot U^1(n)\bigr)$.
\end{enumerate}

Working in $M$, for $p\leq p_0$, $n\in\omega$, and $i<2$, let
\[
 V_{p,n,i}=\{x\in X:(\exists q\leq p)\,
                       q\Vdash\check x\in\dot U^i(\check n)\}.
\]
Each $V_{p,n,i}$ is open in $X$.  Indeed, for $x\in V_{p,n,i}$,
choose $q\leq p$ forcing $\check x\in\dot U^i(\check n)$.
Since $q\Vdash``\dot U^i(\check n) \text{ is open}''$, there exist $q_0\leq q$ and $t\in2^{<\omega}$ with $t\subseteq x$ and
\[
 q_0\Vdash\check{X}\cap[\check t]\subseteq\dot U^i(\check n).
\]
Then $X\cap[t]$ is a neighborhood of $x$ contained in $V_{p,n,i}$.

If $d\in D\setminus\operatorname{dom}(p)$, let
$q=p\cup\{(d,i)\}$.
Then $q\leq p$ and $q\Vdash \dot g(\check d)=\check i$.  Hence $q\Vdash\check d\in\dot U^i(\check n)$,
so
\[
 D\setminus\operatorname{dom}(p)\subseteq V_{p,n,i}.
\]
For $x\in X\setminus D$, by \ref{item:intersection} there exist $p\leq p_0$, $n\in\omega$, and $i<2$ such that $p\Vdash \check x\notin\dot U^i(\check n)$.
So
$x\notin V_{p,n,i}$.  Also $x\notin\operatorname{dom}(p)$,
since $\operatorname{dom}(p)\subseteq D$.

Working in $M$, for each $p\leq p_0$ and $m\in\omega$, let
\[
 W_{p,m}=X\cap\bigcup_{d\in\operatorname{dom}(p)}[d\restriction m].
\]
These sets are clopen in $X$ and, since $\operatorname{dom}(p)$ is
finite, $\bigcap_{m\in\omega}W_{p,m}=\operatorname{dom}(p)$.
The preceding observations give, in $M$,
\[
 D=\bigcap\bigl\{V_{p,n,i}\cup W_{p,m}:p\leq p_0,\ n,m\in\omega,\ i<2\bigr\}.
\]
Since $\Fn(D,2)$ is countable in $M$, this is a countable
intersection of open subsets of $X$ in $M$.
Thus $D$ is $G_\delta$ in $X$ in $M$, a contradiction.
\end{proof}

\begin{corollary}\label{cohen:equivalence}
      Let $X\subseteq\Cantor$.
      \begin{enumerate}[label=(\roman*),ref=(\roman*)]
            \item\label{cohen:equivalence:preserve} If $X$ is a $\lambda$-set, then $\mathbb C\Vdash\text{$\check X$ is a $\lambda$-set}$.
            \item\label{cohen:equivalence:destroy} If $X$ is not a $\lambda$-set, then $\mathbb C\Vdash\text{$\check X$ is not a weak $\lambda$-set}$.
      \end{enumerate}
\end{corollary}

\begin{proof}
For \ref{cohen:equivalence:preserve}, we must see that for every $p\in\mathbb C$ and $\mathbb C$-name $\dot D$ such that $p\Vdash\bigl(\dot D\subseteq\check X\text{ and }|\dot D|\leq\aleph_0\bigr)$, we have $p\Vdash``\dot D$ is relatively $G_\delta$ in $\check X$''.

By Lemma~\ref{cohen:countable-cover}, there exists a countable $A\subseteq X$ such that $p\Vdash\dot D\subseteq\check A$.  Since $X$ is a $\lambda$-set, $A$ is relatively $G_\delta$ in $X$, and hence $p\Vdash``\check A$ is relatively $G_\delta$ in $\check X$''.
As also $p\Vdash|\check A\setminus \dot D|\leq\aleph_0$, it follows that $p\Vdash``\dot D$ is relatively $G_\delta$ in $\check X$''.

For \ref{cohen:equivalence:destroy}, let $D\subseteq X$ be countable and not relatively $G_\delta$ in $X$.  By Lemma~\ref{cohen:coloring},
\[
      \mathbb C\Vdash
      \begin{gathered}
      \exists\,B\subseteq\check D\,\left(\text{no relatively }
      F_\sigma\text{ and }G_\delta\text{ subset }H\text{ of }\check X
      \text{ satisfies }H\cap\check D=B\right).
      \end{gathered}
\]
Thus, $\mathbb C\Vdash\text{$\check X$ is not a weak $\lambda$-set}$.
\end{proof}

A \emph{Luzin set} is an uncountable $L\subseteq\Cantor$ whose
intersection with every meager subset of $\Cantor$ is countable.
By \cite[Theorem~4.4]{RodriguesRonchimSzeptycki}, every weak
$\lambda$-set contained in $\Cantor$ is perfectly meager, hence
meager in $\Cantor$.
A Luzin set $L$ is not meager: otherwise the Luzin property applied
to $L$ itself would make $L$ countable.  Thus no Luzin set is a weak
$\lambda$-set.

The following result is due to Miller.
It is also cited in the proof of
Theorem~22 of Miller's survey~\cite{MillerSurvey}.

\begin{proposition}[Miller~{\cite[pp.~577--578]{Miller1983}}]
\label{prop:cohen-miller}
Assume $\CH$.  Then:
\[
 \mathbb C_{\omega_2}\Vdash
 \text{every }X\subseteq\Cantor\text{ of size }\aleph_2
 \text{ contains a continuous injective image of a Luzin set}.
\]
\end{proposition}

In Proposition~\ref{prop:cohen-miller}, both the Luzin set and the set
containing its continuous injective image are subsets of $\Cantor$.
In the following lemma, a weak $\lambda$-set may be any
separable metrizable space.

\begin{lemma}
\label{cohen:luzin-bound}
Assume that every $X\subseteq\Cantor$ of size $\aleph_2$ contains
a continuous injective image of a Luzin set.  Then every weak
$\lambda$-set has size at most $\aleph_1$.
\end{lemma}

\begin{proof}
Suppose that a weak $\lambda$-set $X$ has size
at least $\aleph_2$.  Identify $X$ with a subspace of $[0,1]^\omega$
and fix a continuous surjection $\pi:\Cantor\to[0,1]^\omega$.
Choose $Y\subseteq\Cantor$ of size $\aleph_2$ such that
$\pi[Y]\subseteq X$ and $\pi\restriction Y$ is injective.
By hypothesis, there is a continuous injection
$f:L\to Y$ from a Luzin set $L$.  Thus $g=\pi\circ f:L\to X$ is
continuous and injective.  Pulling back relatively $F_\sigma$ and
$G_\delta$ separators along $g$ shows that $L$ is a weak
$\lambda$-set, contrary to the fact that no Luzin set is a weak
$\lambda$-set.
\end{proof}

\begin{corollary}
\label{cohen:small}
Assume $\CH$.  Then:
\[
 \mathbb C_{\omega_2}\Vdash
 \text{every weak $\lambda$-set has size at most $\aleph_1$}.
\]
\end{corollary}

\begin{proof}
Apply Lemma~\ref{cohen:luzin-bound} using Proposition~\ref{prop:cohen-miller}.
\end{proof}

The following lemma is well-known and follows from standard techniques involving the use of nice names;
see \cite[Lemma~IV.3.10 and the proof of Lemma~V.2.19]{Kunen}.

\begin{lemma}
\label{cohen:support}
Let $p\in\mathbb C_{\omega_2}$ and let $\dot X$ be a name such that
\[
 p\Vdash\varnothing\ne\dot X\subseteq\Cantor
              \text{ and }|\dot X|\leq\aleph_1.
\]
There exist $S\subseteq\omega_2$ of size at most $\aleph_1$, containing
all coordinates of $p$, and a $\Fn(S\times\omega,2)$-name $\dot Y$
such that $p\Vdash\dot X=\dot Y$.
\end{lemma}

We are now ready to present a proof of Theorem~\ref{thm:intro-cohen-weak-lambda}.

\begin{theoremproof}{thm:intro-cohen-weak-lambda}
Work over a countable transitive model $M$ of $\ZFC+\CH$.
Put
$\kappa=(\omega_2)^M$, and let $G$ be $\mathbb C_\kappa$-generic
over $M$.  Working in $M[G]$, fix a nonempty weak $\lambda$-set
$X$.  By Corollary~\ref{cohen:small}, $|X|\leq\aleph_1<\mathfrak c$.
Thus $X$ is zero-dimensional, as observed in
Section~\ref{sec:weak-lambda}.  Since $X$ is separable and
metrizable, we may assume that $X\subseteq\Cantor$.

Choose a $\mathbb C_\kappa$-name $\dot X\in M$ with $\dot X_G=X$ and $p\in G$ such that
$p\Vdash\varnothing\ne\dot X\subseteq\Cantor$ and
$p\Vdash|\dot X|\leq\aleph_1$.
Working in $M$, apply Lemma~\ref{cohen:support} to obtain
$S\subseteq\kappa$ of size at most $\aleph_1$ and a
$\Fn(S\times\omega,2)$-name $\dot Y$ such that
$p\Vdash_{\mathbb C_\kappa}\dot X=\dot Y$.
Choose $\beta\in\kappa\setminus S$ and put
\[
 P=\Fn((\kappa\setminus\{\beta\})\times\omega,2),\qquad
 Q=\Fn(\{\beta\}\times\omega,2).
\]
Let $H=G\cap P$ and $K=G\cap Q$.  By
\cite[Lemma~V.2.18]{Kunen}, $H$ is $P$-generic over $M$, $K$ is
$Q$-generic over $M[H]$, and $M[G]=M[H][K]$.
Put $N=M[H]$.  Since $\dot Y$ is a $P$-name and $p\in G$,
$X=\dot Y_G=\dot Y_H\in N$.

The poset $Q$ is naturally isomorphic to $\mathbb C$.
Assume, towards a contradiction, that $N\models\text{$X$ is not a $\lambda$-set}$.
Corollary~\ref{cohen:equivalence}\ref{cohen:equivalence:destroy},
applied in $N$, implies that $X$ is not a weak $\lambda$-set in
$N[K]=M[G]$, a contradiction.
Hence $N\models\text{$X$ is a $\lambda$-set}$, and
Corollary~\ref{cohen:equivalence}\ref{cohen:equivalence:preserve},
again applied in $N$, implies that $X$ is a $\lambda$-set in
$M[G]$.
\end{theoremproof}

\begin{corollary}
\label{cohen:independence}
Assuming $\operatorname{Con}(\ZFC)$, the existence of a weak
$\lambda$-set which is not a $\lambda$-set is independent of
$\ZFC+\mathfrak b=\mathfrak s=\aleph_1$.
\end{corollary}

\begin{proof}
Working in the metatheory, it is known that $\operatorname{Con}(\ZFC)$ implies $\operatorname{Con}(\ZFC+\CH)$ since $\CH$ holds in the constructible universe $\mathbf L$ (see e.g.~\cite[Theorem~II.6.24]{Kunen}).
Since $\CH$ implies $\mathfrak p=\mathfrak c$ and
$\mathfrak b=\mathfrak s=\aleph_1$, under $\CH$
Theorem~\ref{thm:intro-p-c-weak-lambda} gives a weak $\lambda$-set
which is not a $\lambda$-set.
So $\ZFC+\mathfrak b=\mathfrak s=\aleph_1$ does not imply that every weak $\lambda$-set is a $\lambda$-set.

On the other hand, working in $\ZFC+\CH$, Theorem~\ref{thm:intro-cohen-weak-lambda} shows that
\[
      \mathbb C_{\omega_2}\Vdash
      \text{every weak $\lambda$-set is a $\lambda$-set}.
\]
As $\mathbb C_{\omega_2}$ also forces $\mathfrak b=\mathfrak s=\aleph_1$
(see e.g.~\cite[Section~11.3]{Blass}), $\Con(\ZFC+\CH)$ implies
the consistency of $\ZFC+\mathfrak b=\mathfrak s=\aleph_1$ together with
the assertion that every weak $\lambda$-set is a $\lambda$-set.
Thus, assuming $\Con(\ZFC)$, the existence of a weak $\lambda$-set which is not a $\lambda$-set cannot be proved from $\ZFC+\mathfrak b=\mathfrak s=\aleph_1$.
\end{proof}

\section{An almost-normal MAD family under CH}
\label{sec:almost-normal-mad-ch}

In this section we construct an almost-normal MAD family under the Continuum Hypothesis.
This construction improves a previous result appearing in \cite{CorralPublished}, where the existence of an almost-normal almost disjoint family which fails to be normal is shown to follow from {\sf CH}. An earlier unpublished version of that paper (see \cite{CorralPreprint}) claimed that it actually follows from {\sf CH} that almost-normal MAD families exist.
However, there was a gap in the argument.
On the other hand, the main idea of the construction there was the use of good families of partitioners to obtain maximality.
We now fix that gap and provide a correct proof of the existence of an almost-normal MAD family under {\sf CH} using good families of partitioners.

\vspace{3mm}

We start with some notation, terminology and definitions. For each set $D\subseteq\omega$ let $D^0=\omega\setminus D$ and $D^1=D$.
For any collection $\mathcal H\subseteq\Pow(\omega)$ and any set
$C\subseteq\omega$, we write
\[
 \mathcal H_C=\{a\in\mathcal H:|a\cap C|=\omega\}.
\]

For each $s\in\Fn(\mathcal H,2)$, define the Boolean combination determined by $s$ as
\[
\mathcal H^s:=\bigcap_{H\in\operatorname{dom}(s)}H^{s(H)}.
\]
We define $\mathcal H^\varnothing=\omega$.
Finally, if $\mathcal A,\cD\subseteq\mathcal P(\omega)$ and $s\in\Fn(\cD,2)$, we let
$$\mathcal A(s)=\{a\in\mathcal A:a\subseteq^*\cD^s\}.$$

The observation below follows directly from the definitions.
Its proof is straightforward and left to the reader.

\begin{remark}
\label{obs:cell-properties}
Let $\D\subseteq\mathcal P(\omega)$, $\mathcal H\subseteq[\omega]^\omega$ and $s, t\in\Fn(\D,2)$.

\begin{enumerate}[label=(\alph*),ref=(\alph*)]
    \item If $s\subseteq t$, then $\D^t\subseteq\D^s$ and $\mathcal H(t)\subseteq\mathcal H(s)$.
    \item If $s\perp t$ (that is, if $s$ and $t$ are incompatible as functions), then $\D^s\cap\D^t=\varnothing$ and $\mathcal H(s)\perp\mathcal H(t)$.
    \item If $s, t$ are compatible, then $\D^{s\cup t}=\D^s\cap\D^t$ and $\mathcal H(s\cup t)=\mathcal H(s)\cap\mathcal H(t)$.
    \item For $D\in \D$, $i\in 2$, and $s\in\Fn(\D,2)$ with $D\notin\operatorname{dom}(s)$,
    $\D^{s\cup\{(D,i)\}}=\D^s\cap D^{i}$.
    \item If $\mathcal D$ is a collection of partitioners for an almost disjoint family $\mathcal A$, then $\cD^s$ is a partitioner for each $s\in\Fn(\mathcal D,2)$.
\end{enumerate}
\end{remark}

\begin{definition}
Let $\D$ be a collection of partitioners for an almost disjoint family $\A$.
We say that $\D$ is
\emph{good for $\A$} if, for every $s\in\Fn(\D,2)$,
\[
 |\A(s)|<\omega
 \quad\Longrightarrow\quad
 \D^s\eqstar\bigcup\A(s).
\]
\end{definition}

Note that if $\D$ is good for $\A$, then every subfamily of $\D$ is also good for $\A$.
The next lemma gives a condition ensuring that a partitioner for a given partition of $\A$ preserves goodness when added to a good family of partitioners $\D$.

\begin{lemma}
\label{lem:goodness-separation-constraints}
Let $\A$ be an almost disjoint family, $\D$ be a good family of partitioners
for $\A$, and partition $\A=\B\mathbin{\dot\cup}\mathcal E$.
For each $s\in\Fn(\D,2)$, let
\[
 L_s=\D^s\setminus\bigcup\mathcal E(s),
 \qquad
 R_s=\D^s\setminus\bigcup\B(s)
\]
and set
\[
\begin{aligned}
 \mathscr L&=\B\cup
 \{L_s:s\in\Fn(\D,2),\ |\mathcal E(s)|<\omega\},\\
 \mathscr R&=\mathcal E\cup
 \{R_s:s\in\Fn(\D,2),\ |\B(s)|<\omega\}.
\end{aligned}
\]
Then $\mathscr L\perp\mathscr R$ and every $D\subseteq\omega$ that separates $\mathscr L$ from $\mathscr R$ also separates $\B$ from $\mathcal E$ and makes $\D\cup\{D\}$ good for $\A$.
\end{lemma}

\begin{proof}
Since $\B$ and $\mathcal E$ are disjoint subfamilies of $\A$,
$\B\perp\mathcal E$.  Fix $s\in\Fn(\D,2)$ with
$|\B(s)|<\omega$ and let $a\in\B$. Then either $a\cap \cD^s$ is finite or $a\in\cB(s)$, in both cases $|a\cap R_s|<\omega$.
Similarly, for every $s\in\Fn(\D,2)$ with $|\mathcal E(s)|<\omega$
and every $e\in\mathcal E$, we have $|L_s\cap e|<\omega$.

It remains to consider $L_s\cap R_t$.  Fix $s,t\in\Fn(\D,2)$ such that
$\mathcal E(s)$ and $\B(t)$ are finite.
If $s,t$ are incompatible, then
$L_s\cap R_t\subseteq\D^s\cap\D^t=\varnothing$.
If they are compatible, then
\[
 \A(s\cup t)
 =\mathcal E(s\cup t)\cup\B(s\cup t)
 \subseteq\mathcal E(s)\cup\B(t).
\]
The right-hand side is finite, so goodness gives
\[
 \cD^{s\cup t}=\D^s\cap\D^t
 \subseteq^*\bigcup\mathcal E(s)\cup\bigcup\B(t).
\]
Hence $L_s\cap R_t\subseteq^*\bigcup\mathcal E(s)\cup\bigcup\B(t)$ which from their definition implies that $L_s\cap R_t$ is finite.

\vspace{3mm}

Let us now show that $\cD\cup\{D\}$ is a good family of partitioners for each $D$ that separates $\mathscr L$ from $\mathscr R$.
For such a set, we also have that $D$ separates $\B$ from $\mathcal E$, so it is a partitioner for $\A$. It remains to show that goodness is preserved when adding $D$.
If $D\in\D$, then $\D\cup\{D\}=\D$ and we are done.  Assume $D\notin\D$ and put
$\D'=\D\cup\{D\}$.  Fix $s\in\Fn(\D,2)$ and put
$s_i=s\cup\{(D,i)\}$ for $i\in2$.
Since $D$ separates $\B$ from $\mathcal E$, we have
\[
 \A(s_1)=\B(s)
 \qquad\text{and}\qquad
 \A(s_0)=\mathcal E(s).
\]
If $\B(s)$ is finite, then clearly $\bigcup\cB(s)\subseteq^* D\cap\cD^{s}=\D'^{s_1}$, as each $b\in\cB(s)$ is almost contained in $\D'^{s_1}$ and $\cB(s)$ is finite. On the other hand, since $R_s\in\mathscr R$, we have that $R_s\cap D=\big(\cD^s\setminus\bigcup\cB(s)\big)\cap D$ is finite, or equivalently, that $\D'^{s_1}=\cD^s\cap D\subseteq^*\bigcup \cB(s)$. Together, these inclusions give
\[
 \D'^{s_1}\eqstar\bigcup\B(s)=\bigcup\cA(s_1).
\]
If $\mathcal E(s)$ is finite, a dual argument gives that
$\D'^{s_0}\eqstar\bigcup\cA(s_0)$.
We conclude that $\D\cup\{D\}$ is good for $\A$.
\end{proof}

The next lemma gives a sufficient condition for adding a separator
containing a prescribed set $C\subseteq\omega$ while preserving goodness.

\begin{lemma}
\label{lem:add-separator}
Let $\A$ be countable and almost disjoint, let $\D$ be a countable good
family of partitioners for $\A$, and let $C\subseteq\omega$.
Suppose that, for every $s\in\Fn(\D,2)$,
\begin{equation}
\label{eq:separator-criterion}
 |\A_C(s)|<\omega
 \quad\Longrightarrow\quad
 C\cap\D^s\subseteq^*\bigcup\A_C(s).
\end{equation}
Then there exists a set $D\subseteq\omega$ with $C\subseteq D$ separating
$\A_C$ from $\A\setminus\A_C$ such that $\D\cup\{D\}$ is good for $\A$.
\end{lemma}

\begin{proof}
Let $\mathscr L,\mathscr R$ be the collections from
Lemma~\ref{lem:goodness-separation-constraints} with $\B=\A_C$ and
$\mathcal E=\A\setminus\A_C$.  They are orthogonal and, since
$\A$ and $\Fn(\D,2)$ are countable, both collections are countable.
By the definition of $\A_C$, $C$ meets every member of $\mathcal E$
only finitely.  For each $s\in\Fn(\D,2)$ with
$|\A_C(s)|<\omega$, \eqref{eq:separator-criterion} gives
\[
 \left|C\cap\left(\D^s\setminus\bigcup\A_C(s)\right)\right|<\omega.
\]
Thus $\{C\}\perp\mathscr R$, and hence
$(\{C\}\cup\mathscr L)\perp\mathscr R$.
Lemma~\ref{lem:countable-separation} gives a set $D\supseteq C$ separating
$\{C\}\cup\mathscr L$ from $\mathscr R$.
By Lemma \ref{lem:goodness-separation-constraints},
$D$ separates $\A_C$ from its complement in $\A$ and
$\D\cup\{D\}$ is good for $\A$.
\end{proof}

\begin{lemma}
\label{lem:goodness-preserved}
Suppose that $\D$ is good for $\A$, that $\A\subseteq\A'$, and that
$\A'$ is almost disjoint.  If every $D\in\D$ is a partitioner for $\A'$,
then $\D$ is good for $\A'$.
\end{lemma}

\begin{proof}
Fix $s\in\Fn(\D,2)$ and suppose that $\A'(s)$ is finite.
Then $\A(s)$ is finite, and goodness for $\A$ gives
$\D^s\eqstar\bigcup\A(s)$. Almost disjointness of $\cA'$ implies that $b\cap\cD^s$ is finite for each $b\in\A'\setminus\A$.
Thus $\A'(s)=\A(s)$ and $\D^s\eqstar\bigcup\A'(s)$.
\end{proof}

We now include some terminology in order to simplify our coming recursive construction.

\begin{definition}
Let $\mathscr P\subseteq \Pow(\omega)\times\Pow(\omega)$ be countable.
We say that $(\cA,\mathscr P)$ is a \emph{countable compatible system} if:
\begin{enumerate}[label=(\alph*)]
\item $\cA$ is a countable almost disjoint family,
\item $\D=\{D:(C,D)\in\mathscr P\}$ is a good family of partitioners for $\A$;
\item $C\subseteq^*D$ for every $(C,D)\in\mathscr P$ and
\item\label{item:compatibility} for each $(C,D)\in\mathscr P$,
$\cA_C=\{a\in\cA:a\subseteq^* D\}$.
\end{enumerate}
\end{definition}

\begin{lemma}
\label{lem:one-point-extension}
Let $(\A,\mathscr P)$ be a countable compatible system, and let
$X\in[\omega]^\omega$ be almost disjoint from $\A$.  There exists
$a_*\in[\omega]^\omega$ such that $|a_*\cap X|=\omega$ and
$(\A\cup\{a_*\},\mathscr P)$ is a countable compatible system.
\end{lemma}

\begin{proof}
Adding the pair $(\omega,\omega)$ to $\mathscr P$ does not affect compatibility, so we may assume that $\mathscr P\neq\emptyset$.
Enumerate $\mathscr P$ as $\{(C_n,D_n):n<\omega\}$,
allowing repetitions if necessary, and put $\D=\{D_n:n<\omega\}$.

The hypotheses imply that $\A$ is infinite: otherwise goodness applied
to $s=\varnothing$ would make $\omega$ almost covered by finitely many
members of $\A$, contrary to the choice of $X$.  Enumerate
$\A=\{a_j:j<\omega\}$.
Starting with $s_0=\varnothing$, construct an increasing sequence
$(s_n)_{n<\omega}$ in $\Fn(\D,2)$ such that, for each $n \in \omega$,
\begin{equation}
\label{eq:X-cell-infinite}
 \operatorname{dom}(s_n)=\{D_i:i<n\}
 \qquad\text{and}\qquad |X\cap\cD^{s_n}|=\omega.
\end{equation}

This can be easily done by picking $i\in 2$ such that $(X\cap\cD^{s_n})\cap D^{i}_{n}$ is infinite and by defining $s_{n+1}(D_n)=i$.
Define:
\[
 h=\bigcup_{n<\omega}s_n:\mathcal D\to 2
\]

For every $n<\omega$, the set $\A(s_n)$ is infinite, otherwise
goodness would give
 $\D^{s_n}\eqstar\bigcup\A(s_n)$,
where the right-hand side is a finite union of members of $\A$,
contradicting
\eqref{eq:X-cell-infinite} and the fact that $X$ is almost disjoint from
$\A$.

Choose recursively finite nonempty sets $F_n$ for $n<\omega$ such that
\begin{enumerate}[label=(\roman*)]
\item\label{item:Fn-position} $F_n\subseteq\cD^{s_n}\setminus
\left(\bigcup_{j<n}a_j\cup\bigcup_{m<n}F_m\right)$;
\item $F_n\cap X\neq\varnothing$;
\item\label{item:Fn-requirements} for every $i<n$ with $h(D_i)=1$, $F_n\cap C_i\neq\varnothing$.
\end{enumerate}
By \eqref{eq:X-cell-infinite},
$X\cap\cD^{s_n}$ is infinite.
As $X$ is almost disjoint from $\cA$, we can choose $l\in X$ satisfying
item~\ref{item:Fn-position}. Similarly, if
$i<n$ and $h(D_i)=1$, choose
$b_n\in\A(s_n)\setminus\{a_j:j<n\}$.

As $b_n\in\cA(s_n)$ and $s_n(D_i)=h(D_i)=1$, we have that $b_n\subseteq^*D_i$, so
compatibility ensures that
$|b_n\cap C_i\cap\cD^{s_n}|=\omega$. Thus we can pick $l\in (b_n\cap C_i\cap\cD^{s_n})\setminus\bigcup_{j<n}a_j$ out of
$\bigcup_{m<n}F_m$. With these witnesses we can build $F_{n}$.

Let
$a_*=\bigcup_{n<\omega}F_n$.
The $F_n$ are pairwise disjoint and each meets $X$, so
$|a_*\cap X|=\omega$.  For a fixed $j<\omega$, every $F_n$ with $n>j$ avoids
$a_j$, and therefore $|a_*\cap a_j|<\omega$.  Thus
$\A\cup\{a_*\}$ is almost disjoint.

Fix $i<\omega$.  For every $n>i$, $F_n\subseteq D_i^{h(D_i)}$.
If $h(D_i)=1$, then
$a_*\subseteq^*D_i$, and item~\ref{item:Fn-requirements} gives infinitely many
distinct points of $a_*\cap C_i$.  If $h(D_i)=0$, then
$|a_*\cap D_i|<\omega$; since $C_i\subseteq^*D_i$, also
$|a_*\cap C_i|<\omega$.  Thus every $D_i$ remains a partitioner for
$\A\cup\{a_*\}$ and item~\ref{item:compatibility} holds for
$a_*$.  Lemma~\ref{lem:goodness-preserved} gives goodness for the extension,
so the resulting system is compatible.
\end{proof}

We now extend the almost disjoint family to force it to satisfy
\eqref{eq:separator-criterion} for a prescribed set $C\subseteq\omega$,
preserving compatibility with all pairs in a given countable compatible
system. 

\begin{lemma}
\label{lem:saturation}
Let $(\A,\mathscr P)$ be a countable compatible system, and let
$C\subseteq\omega$.  Put $\D=\{D:(C',D)\in\mathscr P\}$.
There exists a countable almost disjoint family $\mathcal A'$ extending $\A$ such
that:
\begin{enumerate}[label=(\alph*),ref=(\alph*)]
\item $(\A',\mathscr P)$ is a countable compatible system;
\item\label{item:saturation-conclusion} for every $s\in\Fn(\D,2)$,
if $|\A'_C(s)|<\omega$, then $C\cap\cD^s\subseteq^*\bigcup\A'_C(s)$.
\end{enumerate}
\end{lemma}

\begin{proof}
Enumerate $\Fn(\D,2)$ as $(s_k)_{k<\omega}$ in which every
function occurs infinitely often. We construct an increasing sequence of
countable almost disjoint families
\[
 \A_0=\A\subseteq\A_1\subseteq\cdots
\]
such that $(\A_k,\mathscr P)$ is a countable compatible system
for every $k<\omega$.

At stage $k<\omega$, suppose that $\A_k$ has been constructed.  Put
\[
 \B_k=(\A_k)_C(s_k)
 \qquad\textnormal{ and }\qquad
 X_k=(C\cap\cD^{s_k})\setminus\bigcup\B_k.
\]
If $\B_k$ is infinite or $X_k$ is finite, set
$\A_{k+1}=\A_k$.  Otherwise, $\B_k$ is finite and $X_k$ is infinite.
We have that $X_k$ is almost disjoint from $\A_k$ and so we can apply Lemma~\ref{lem:one-point-extension} to
$(\A_k,\mathscr P)$ and $X_k$.  Let $b_k$ be the resulting
member and set $\A_{k+1}=\A_k\cup\{b_k\}$.
The extended system is compatible and
$|b_k\cap X_k|=\omega$ with $X_k\subseteq C\cap\cD^{s_k}$.
Since every member of $\D$ is a partitioner for $\A_{k+1}$ and $b_k\in\cA_{k+1}$, we necessarily have that $b_k\subseteq^*\cD^{s_k}$.
Thus $b_k\in(\A_{k+1})_C(s_k)$.  Also $b_k\notin\A_k$, since $X_k$ is
almost disjoint from $\A_k$.

Let $\A'=\bigcup_{k<\omega}\A_k$.
This set is countable and almost disjoint.  Every $D\in\D$ is a
partitioner for $\A'$ and the condition in item~\ref{item:compatibility} holds for every
member of $\A'$, since each member belongs to some $\A_k$.
Lemma~\ref{lem:goodness-preserved}, applied to
$\A\subseteq\A'$, shows that $\D$ is good for $\A'$.
Hence $(\A',\mathscr P)$ is a countable compatible system.

Fix $s\in\Fn(\D,2)$ and suppose that $\A'_C(s)$ is finite.
Choose $k$ with $s_k=s$ and $\A'_C(s)\subseteq\A_k$.
Then $\B_k=\A'_C(s)$.  If $X_k$ were infinite, the construction would
add $b_k\in\A'_C(s)\setminus\A_k$, a contradiction.
Thus $X_k=(C\cap\cD^{s})\setminus\bigcup\A'_C(s)$ is finite, proving
item~\ref{item:saturation-conclusion}.
\end{proof}

Combining the previous results gives the following extension property.

\begin{proposition}
\label{prop:strong-extension}
Let $(\A,\mathscr P)$ be a countable compatible system and
$C\subseteq\omega$.  There exist a countable almost disjoint extension
$\A\subseteq\A'$ and a set $D\subseteq\omega$ with $C\subseteq D$ such
that $(\A',\mathscr P\cup\{(C,D)\})$ is a countable compatible system.
\end{proposition}

\begin{proof}
Put $\D=\{D':(C',D')\in\mathscr P\}$.
Let $\A'$ be the extension given by Lemma~\ref{lem:saturation}. By item~\ref{item:saturation-conclusion} and
Lemma~\ref{lem:add-separator}, there is a set $D\subseteq\omega$ with
$C\subseteq D$ separating $\A'_C$ from $\A'\setminus\A'_C$ and making
$\D\cup\{D\}$ good for $\A'$.
\end{proof}

We are now ready to prove Theorem~\ref{thm:intro-ch-almost-normal-mad}.

\begin{theoremproof}{thm:intro-ch-almost-normal-mad}
By $\CH$, fix an enumeration $\mathcal P(\omega)=\{C_\alpha:\alpha<\omega_1\}$.
Let $\A_0=\{a_n:n<\omega\}$ be a partition of $\omega$ into infinite sets.

We recursively construct an increasing chain of countable
almost disjoint families $\A_\alpha$ and sets $D_\alpha\subseteq\omega$ satisfying the following.
\begin{enumerate}[label=(\roman*)]
\item $\D_\alpha=\{D_\xi:\xi<\alpha\}$ is a good family of partitioners for $\A_\alpha$;
\item $C_\xi\subseteq D_\xi$ for every $\xi<\alpha$;
\item\label{item:mad-recursion-compatibility} for each $\xi<\alpha$, $(\cA_\alpha)_{C_\xi}=\{a\in\cA_\alpha:a\subseteq^* D_\xi\}$
\end{enumerate}

Suppose that $\alpha<\omega_1$ and that the construction has been carried
out through stage $\alpha$.
We show how to extend the construction to stage $\alpha+1$.

Apply Proposition~\ref{prop:strong-extension} to $(\A_\alpha,\{(C_\xi,D_\xi):\xi<\alpha\})$ and
$C_\alpha$ to obtain a countable almost disjoint family
$\A_{\alpha+1}$ extending $\A_\alpha$ and a set $D_\alpha\supseteq C_\alpha$
such that $\D_{\alpha+1}=\D_\alpha\cup\{D_\alpha\}$ is good for
$\A_{\alpha+1}$ and, for every $\xi\leq\alpha$, $(\A_{\alpha+1})_{C_\xi}
 =\{a\in\A_{\alpha+1}:a\subseteq^*D_\xi\}$.
 Thus, all three requirements hold
at stage $\alpha+1$.

If $\gamma<\omega_1$ is a limit ordinal, define $\A_\gamma=\bigcup_{\alpha<\gamma}\A_\alpha$.
This family is countable and almost disjoint. For each $\xi<\gamma$, $D_\xi$ is a
partitioner for $\A_\gamma$, since every member of $\A_\gamma$ belongs to some
$\A_\alpha$ with $\xi<\alpha<\gamma$. Moreover,
condition~\ref{item:mad-recursion-compatibility} holds, since
\[
\begin{aligned}
 (\A_\gamma)_{C_\xi}
 &=\{a\in\A_\gamma:|a\cap C_\xi|=\omega\}
 =\bigcup_{\xi<\alpha<\gamma}\{a\in\A_\alpha:|a\cap C_\xi|=\omega\}\\
 &=\bigcup_{\xi<\alpha<\gamma}\{a\in\A_\alpha:a\subseteq^*D_\xi\}
 =\{a\in\A_\gamma:a\subseteq^*D_\xi\}.
\end{aligned}
\]
To see that (i) holds, fix $s\in\Fn(\D_\gamma,2)$.
Choose $\beta<\gamma$ with
$\operatorname{dom}(s)\subseteq\D_\beta$.
Since $\D_\beta$ is good for $\A_\beta$ and every member of $\D_\beta$ is a
partitioner for $\A_\gamma$, Lemma~\ref{lem:goodness-preserved} implies that
$\D_\beta$ is good for $\A_\gamma$. Hence, if $|\A_\gamma(s)|<\omega$, then
$\D_\gamma^s=\D_\beta^s\eqstar\bigcup\A_\gamma(s)$.
This completes the recursive construction.

Set $\A=\bigcup_{\alpha<\omega_1}\A_\alpha$.
We claim that $\A$ is a MAD family.
Indeed, given
$X\in[\omega]^\omega$, choose $\alpha<\omega_1$ such that $X=C_\alpha$,
and put $s=\{(D_\alpha,1)\}\in\Fn(\D_{\alpha+1},2)$.
If no member of $\A_{\alpha+1}$ met $X$ infinitely, then, by compatibility,
$\A_{\alpha+1}(s)=(\A_{\alpha+1})_{C_\alpha}=\varnothing$.
Goodness applied to $s$, with $\cD_{\alpha+1}^{s}=D_\alpha$, would then make
$D_\alpha$ finite, contradicting $X\subseteq D_\alpha$.
Thus some member of $\A_{\alpha+1}$, and hence of $\A$, meets $X$ infinitely.

Finally, fix $C\subseteq\omega$ and choose $\alpha<\omega_1$ with $C=C_\alpha$.
Given $a\in\A$, we can pick $\gamma<\omega_1$ such that $\gamma>\alpha$ and
$a\in\A_\gamma$. 
At stage $\gamma$, $D_\alpha$ is a partitioner and $(\cA_\gamma)_{C_\alpha}=\{a\in\cA_\gamma:a\subseteq^* D_\alpha\}$ gives
\[
 a\subseteq^*D_\alpha
 \quad\Longleftrightarrow\quad
 |a\cap C|=\omega.
\]
Since $a\in\A$ was arbitrary, $D_\alpha$ separates $\A_C$ from
$\A\setminus\A_C$.  Proposition~\ref{prop:almost-normal-criterion}
now implies that $\Psi(\A)$ is almost-normal.
\end{theoremproof}

Theorem~\ref{thm:intro-ch-almost-normal-mad} gives the consistency,
relative to $\ZFC$, of the existence of an almost-normal MAD family. We now turn to its independence.

\section{No almost-normal MAD families under
\texorpdfstring{$\OCA+\mathfrak{ap}>\omega_1$}{OCA + ap > omega1}}
\label{sec:oca-obstruction}

In \cite{CorralPublished}, the author proved that the Proper Forcing Axiom (PFA) implies that no
almost-normal MAD family exists.  We obtain the
same conclusion under $\OCA+\mathfrak{ap}>\omega_1$, where $\OCA$ is
Todor\v cevi\'c's Open Coloring Axiom, as stated in
\cite[Definition~V.6.25]{Kunen}.
The cardinal $\mathfrak{ap}$ is the least size of an almost disjoint
family $\A$ with a subset $\B\subseteq\A$ for which no
$X\subseteq\omega$ meets every member of $\B$ infinitely and every member
of $\A\setminus\B$ only finitely \cite[p.~14]{BrendleDow}.

A partially ordered set $\mathbb P$ is said to be \emph{$\sigma$-centered} if it can be expressed as a countable union of centered sets.
Luzin gaps in the sense used here were introduced by Todor\v cevi\'c
\cite{TodorcevicAnalyticGaps}.
Given an uncountable almost disjoint family $\A$, one defines its splitting
partial order, denoted by $\mathbb P_\A$.
We will not define $\mathbb P_\A$ nor Luzin gaps in this paper as we will not
use these definitions directly. 
The interested reader may consult e.g.
\cite[Definition~5 and 12]{GuzmanHrusakKoszmider}.
We use only the
following statements:
\begin{enumerate}[label=(\roman*)]
\item If $\A$ is MAD, then $\mathbb P_\A$ is not $\sigma$-centered
\cite[Theorem~9(2)(b)]{GuzmanHrusakKoszmider}.
\item Under $\OCA$, $\mathbb P_\A$ is either $\sigma$-centered or fails
the ccc \cite[Theorem~9(3)(b)]{GuzmanHrusakKoszmider}.
\item Under $\OCA$, if $\mathbb P_\A$ fails the ccc, then $\A$ contains
a $2$-Luzin gap \cite[Theorem~9(1)(c)]{GuzmanHrusakKoszmider}.
\item The two sides of a $2$-Luzin gap are disjoint subfamilies of
cardinality $\omega_1$ which cannot be separated
\cite[p.~1365]{HrusakGuzman}.
\end{enumerate}

From this, the following proposition follows immediately.

\begin{proposition}
\label{prop:oca-inseparable-subfamilies}
Assume $\OCA$.  Every MAD family $\A$ contains disjoint subfamilies
$\B,\C\in[\A]^{\aleph_1}$ which cannot be separated.
\end{proposition}

\begin{proof}
If $\A$ is MAD, (i) and (ii) imply that $\mathbb P_\A$ fails the ccc.
Now apply (iii) and (iv).
\end{proof}

We can now prove Theorem~\ref{thm:intro-oca-ap-no-almost-normal}.

\begin{theoremproof}{thm:intro-oca-ap-no-almost-normal}
Suppose that $\A$ is MAD and $\Psi(\A)$ is almost-normal.  Fix
$\B,\C$ as in Proposition~\ref{prop:oca-inseparable-subfamilies}.
Since $|\B\cup\C|=\omega_1<\mathfrak{ap}$, there is $X\subseteq\omega$
which weakly separates $\B$ from $\C$.  Thus $\B\subseteq\A_X$ and
$\C\subseteq\A\setminus\A_X$.
However, by Proposition~\ref{prop:almost-normal-criterion} there exists a partitioner
separating $\A_X$ from $\A\setminus\A_X$, and hence separating $\B$
from $\C$, a contradiction.
\end{theoremproof}

Since $\mathfrak p\leq\mathfrak{ap}$ \cite[p.~14]{BrendleDow}, the
conclusion also holds under $\OCA+\mathfrak p>\omega_1$, and in
particular under $\OCA+\MA_{\omega_1}$.  As $\MA+\OCA+\mathfrak c=\aleph_2$ is consistent with $\ZFC$ \cite{TodorcevicPartitionProblems}, the above observation,
together with the results of Section~\ref{sec:almost-normal-mad-ch}, shows the following.

\begin{corollary}
    The existence of almost-normal MAD families is independent of $\ZFC$.
\end{corollary}

\section*{Acknowledgements}

The first author was financed, in part, by the S\~ao Paulo Research
Foundation (FAPESP), Brasil. Process Number 2025/07302-0.

The second author was financed, in part, by the S\~ao Paulo Research
Foundation (FAPESP), Brasil. Process Number 2025/07302-0.

The third author was financed, in part, by the S\~ao Paulo Research
Foundation (FAPESP), Brasil. Process Number 2025/08683-7.

\section*{Use of artificial intelligence}

The main mathematical results presented in this paper had the first draft of their proofs generated by
OpenAI's GPT-6 Astra.
The authors carefully reviewed and verified the correctness of these proofs, rewriting them as necessary to ensure their accuracy and clarity.
OpenAI's GPT-6.1 Sol and OpenAI's GPT-6.1 Terra aided the authors with subsequent revisions of the proofs, grammar, and LaTeX formatting.

\bibliographystyle{plainurl}
\bibliography{includes/references}

@misc{Blass,
  author = {Blass, Andreas},
  title  = {Combinatorial Cardinal Characteristics of the Continuum},
  year   = {2003},
  note   = {Preprint, November 24, 2003},
  url    = {https://sites.lsa.umich.edu/ablass/wp-content/uploads/sites/1471/2025/10/hbk.pdf}
}

@article{Bell,
  author  = {Bell, Murray G.},
  title   = {On the combinatorial principle {$P(\mathfrak c)$}},
  journal = {Fundamenta Mathematicae},
  volume  = {114},
  number  = {2},
  year    = {1981},
  pages   = {149--157},
  doi     = {10.4064/fm-114-2-149-157}
}

@article{BrendleDow,
  author  = {Brendle, J{\"o}rg},
  title   = {Dow's principle and {$Q$}-sets},
  journal = {Canadian Mathematical Bulletin},
  volume  = {42},
  number  = {1},
  year    = {1999},
  pages   = {13--24},
  doi     = {10.4153/CMB-1999-002-2}
}

@misc{CorralPreprint,
  author        = {Corral, C{\'e}sar},
  title         = {Madness and (weak) normality},
  year          = {2021},
  eprint        = {2101.01829},
  archiveprefix = {arXiv},
  primaryclass  = {math.GN},
  doi           = {10.48550/arXiv.2101.01829},
  note          = {arXiv:2101.01829v3}
}

@article{CorralPublished,
  author  = {Corral, C{\'e}sar},
  title   = {Madness and weak forms of normality},
  journal = {Acta Mathematica Hungarica},
  volume  = {165},
  number  = {2},
  year    = {2021},
  pages   = {291--307},
  doi     = {10.1007/s10474-021-01186-y}
}

@article{GarciaBalanSzeptycki,
  author  = {Garcia-Balan, Sergio A. and Szeptycki, Paul J.},
  title   = {Weak normality properties in {$\Psi$}-spaces},
  journal = {Fundamenta Mathematicae},
  volume  = {258},
  number  = {2},
  year    = {2022},
  pages   = {137--151},
  doi     = {10.4064/fm935-11-2021}
}

@article{GuzmanHrusakKoszmider,
  author  = {Guzm{\'a}n, Osvaldo and Hru{\v{s}}{\'a}k, Michael and
             Koszmider, Piotr},
  title   = {Almost disjoint families and the geometry of nonseparable
             spheres},
  journal = {Journal of Functional Analysis},
  volume  = {285},
  number  = {11},
  year    = {2023},
  pages   = {110149},
  note    = {Paper No. 110149, 49 pp.},
  doi     = {10.1016/j.jfa.2023.110149}
}

@incollection{HernandezHrusakMrowkaIsbell,
  author    = {Hern{\'a}ndez-Hern{\'a}ndez, Fernando and
               Hru{\v{s}}{\'a}k, Michael},
  title     = {Topology of {Mr{\'o}wka--Isbell} spaces},
  booktitle = {Pseudocompact Topological Spaces: A Survey of Classic and
               New Results with Open Problems},
  editor    = {Hru{\v{s}}{\'a}k, Michael and
               Tamariz-Mascar{\'u}a, {\'A}ngel and Tkachenko, Mikhail},
  series    = {Developments in Mathematics},
  volume    = {55},
  publisher = {Springer},
  address   = {Cham},
  year      = {2018},
  pages     = {253--289},
  doi       = {10.1007/978-3-319-91680-4_8}
}

@incollection{HrusakAlmostDisjointTopology,
  author    = {Hru{\v{s}}{\'a}k, Michael},
  title     = {Almost disjoint families and topology},
  booktitle = {Recent Progress in General Topology {III}},
  editor    = {Hart, Klaas Pieter and van Mill, Jan and Simon, Petr},
  publisher = {Atlantis Press},
  address   = {Paris},
  year      = {2014},
  pages     = {601--638},
  doi       = {10.2991/978-94-6239-024-9_14}
}

@article{HrusakGuzman,
  author  = {Hru{\v{s}}{\'a}k, Michael and Guzm{\'a}n, Osvaldo},
  title   = {{$n$}-Luzin gaps},
  journal = {Topology and its Applications},
  volume  = {160},
  number  = {12},
  year    = {2013},
  pages   = {1364--1374},
  doi     = {10.1016/j.topol.2013.05.015}
}

@article{Hindman,
  author  = {Hindman, Neil},
  title   = {Finite sums from sequences within cells of a partition of
             {$\mathbb N$}},
  journal = {Journal of Combinatorial Theory, Series A},
  volume  = {17},
  number  = {1},
  year    = {1974},
  pages   = {1--11},
  doi     = {10.1016/0097-3165(74)90023-5}
}

@article{KamburelisWeglorz,
  author  = {Kamburelis, Anastasis and W{\k{e}}glorz, Bogdan},
  title   = {Splittings},
  journal = {Archive for Mathematical Logic},
  volume  = {35},
  number  = {4},
  year    = {1996},
  pages   = {263--277},
  doi     = {10.1007/s001530050044}
}

@article{Miller1983,
  author  = {Miller, Arnold W.},
  title   = {Mapping a set of reals onto the reals},
  journal = {The Journal of Symbolic Logic},
  volume  = {48},
  number  = {3},
  year    = {1983},
  pages   = {575--584},
  doi     = {10.2307/2273449},
  url     = {https://people.math.wisc.edu/~awmille1/res/map.pdf}
}

@article{MillerLambdaSets,
  author  = {Miller, Arnold W.},
  title   = {On {$\lambda'$}-sets},
  journal = {Topology Proceedings},
  volume  = {28},
  number  = {1},
  year    = {2004},
  pages   = {179--187},
  url     = {https://topology.nipissingu.ca/tp/reprints/v28/tp28112.pdf}
}

@incollection{MillerSurvey,
  author    = {Miller, Arnold W.},
  title     = {Special sets of reals},
  booktitle = {Set Theory of the Reals},
  editor    = {Judah, Haim},
  series    = {Israel Mathematical Conference Proceedings},
  volume    = {6},
  publisher = {American Mathematical Society},
  year      = {1993},
  pages     = {415--432},
  url       = {https://people.math.wisc.edu/~awmille1/res/survey.pdf}
}

@article{RodriguesRonchimAlmostNormal,
  author  = {Rodrigues, Vinicius de Oliveira and
             Ronchim, Victor dos Santos},
  title   = {Almost-normality of {Isbell--Mr{\'o}wka} spaces},
  journal = {Topology and its Applications},
  volume  = {288},
  year    = {2021},
  pages   = {107470},
  note    = {Paper No. 107470, 13 pp.},
  doi     = {10.1016/j.topol.2020.107470}
}

@article{RodriguesRonchimSzeptycki,
  author  = {Rodrigues, Vinicius de Oliveira and
             Ronchim, Victor dos Santos and Szeptycki, Paul J.},
  title   = {Special sets of reals and weak forms of normality on
             {Isbell--Mr{\'o}wka} spaces},
  journal = {Commentationes Mathematicae Universitatis Carolinae},
  volume  = {64},
  number  = {1},
  year    = {2023},
  pages   = {109--126},
  doi     = {10.14712/1213-7243.2023.014}
}

@article{Rothberger,
  author  = {Rothberger, Fritz},
  title   = {Sur un ensemble toujours de premi{\`e}re cat{\'e}gorie
             qui est d{\'e}pourvu de la propri{\'e}t{\'e}
             {$\lambda$}},
  journal = {Fundamenta Mathematicae},
  volume  = {32},
  number  = {1},
  year    = {1939},
  pages   = {294--300},
  doi     = {10.4064/fm-32-1-294-300}
}

@article{TallNormalMoore,
  author  = {Tall, Franklin D.},
  title   = {Set-theoretic consistency results and topological theorems
             concerning the normal {Moore} space conjecture and related problems},
  journal = {Dissertationes Mathematicae},
  volume  = {148},
  year    = {1977},
  pages   = {1--53}
}

@article{TodorcevicAnalyticGaps,
  author  = {Todor{\v{c}}evi{\'c}, Stevo},
  title   = {Analytic gaps},
  journal = {Fundamenta Mathematicae},
  volume  = {150},
  number  = {1},
  year    = {1996},
  pages   = {55--66},
  doi     = {10.4064/fm-150-1-55-66}
}

@book{TodorcevicPartitionProblems,
  author    = {Todor{\v{c}}evi{\'c}, Stevo},
  title     = {Partition Problems in Topology},
  series    = {Contemporary Mathematics},
  volume    = {84},
  publisher = {American Mathematical Society},
  address   = {Providence, RI},
  year      = {1989},
  isbn      = {0-8218-5091-1},
  doi       = {10.1090/conm/084}
}

@incollection{vanDouwen,
  author    = {van Douwen, Eric K.},
  title     = {The integers and topology},
  booktitle = {Handbook of Set-Theoretic Topology},
  editor    = {Kunen, Kenneth and Vaughan, Jerry E.},
  publisher = {North-Holland},
  address   = {Amsterdam},
  year      = {1984},
  pages     = {111--167},
  doi       = {10.1016/B978-0-444-86580-9.50006-9}
}

@book{Kunen,
  author    = {Kunen, Kenneth},
  title     = {Set Theory},
  series    = {Studies in Logic},
  volume    = {34},
  publisher = {College Publications},
  address   = {London},
  year      = {2011},
  isbn      = {978-1-84890-050-9}
}

@incollection{Goldstern,
  author    = {Goldstern, Martin},
  title     = {Tools for your forcing construction},
  booktitle = {Set Theory of the Reals},
  editor    = {Judah, Haim},
  series    = {Israel Mathematical Conference Proceedings},
  volume    = {6},
  publisher = {Bar-Ilan University},
  address   = {Ramat Gan},
  year      = {1993},
  pages     = {305--360},
  url       = {https://dmg.tuwien.ac.at/goldstern/www/papers/tools.pdf}
}

@misc{osipovweaklambda,
      title={A weak $\lambda$-set that is not a $\lambda$-set under $\mathfrak{p} = \mathfrak{c}$},
      author={Alexander V. Osipov and Anton E. Lipin},
      year={2026},
      eprint={2611.06312},
      archivePrefix={arXiv},
      primaryClass={math.LO},
      url={https://arxiv.org/abs/2611.06312}, 
}

\end{document}